\documentclass[a4paper]{svjour3}   
\usepackage[margin=1in]{geometry}
\smartqed  
\usepackage{graphicx}
\usepackage{color}

\usepackage{float, enumerate, graphicx, amssymb, amsmath , color, longtable, multirow, lineno,url}

\usepackage[justification=centering]{caption}
\usepackage{amsmath}
\usepackage[active]{srcltx}
\usepackage{pifont}
\usepackage{dsfont}
\usepackage{paralist} 
\usepackage{amssymb}
\usepackage{authblk}
\usepackage{algorithm}
\usepackage{algorithmic}
\usepackage{lineno}
\usepackage{epsfig}
\usepackage{epstopdf}
\usepackage[colorlinks,linkcolor=blue,citecolor=blue,                  bookmarksnumbered=true,
bookmarksopen=true,
colorlinks ]{hyperref}

\usepackage{color}

\numberwithin{theorem}{section}
\numberwithin{lemma}{section}
\numberwithin{corollary}{section}
\numberwithin{proposition}{section}
\numberwithin{definition}{section}
\numberwithin{example}{section}

\numberwithin{rem}{section}
\begin{document}

\title{New subspace minimization conjugate gradient methods based on regularization model for unconstrained optimization
}


\author{Ting Zhao \textsuperscript{1}        \and
        Hongwei Liu \textsuperscript{1}
        \and
        Zexian Liu\textsuperscript{2}
}


\institute{\ding {41}Hongwei Liu   \at
	\email{hwliuxidian@163.com}          
	\and
	Ting Zhao \at
	\email{zhaoting\_0322@163.com}\\
	\and
	Zexian Liu \at
	\email{liuzexian2008@163.com}\\
	\and
	 $^1$  School of Mathematics and  Statistics, Xidian University, Xi'an 710126,  China \\
	 $^2$  State Key Laboratory of Scientific and Engineering Computing, Institute of Computational Mathematics and Scientific/Engineering computing, AMSS, Chinese Academy of Sciences, Beijing, 100190, China.\\
}

\date{Received: date / Accepted: date}

\maketitle

\begin{abstract}
 In this paper, two new subspace minimization conjugate gradient methods based on $p - $regularization models are proposed, where a special scaled norm in $p - $regularization model is analyzed. Different choices for special scaled norm lead to different solutions to the $p - $regularized subproblem. Based on the analyses of the solutions in a two-dimensional subspace, we derive new directions satisfying the sufficient descent condition. With a modified nonmonotone line search, we establish the global convergence of the proposed methods under mild assumptions. $R - $linear convergence of the proposed methods are also analyzed. Numerical results show that, for the CUTEr library, the proposed methods are superior to four conjugate gradient methods, which were proposed by Hager and Zhang (SIAM J Optim 16(1):170-192, 2005), Dai and Kou (SIAM J Optim 23(1):296-320, 2013), Liu and Liu (J Optim Theory Appl 180(3):879-906, 2019) and Li et al. (Comput Appl Math 38(1): 2019), respectively.
\end{abstract}

\section{Introduction}
\label{intro}
Conjugate gradient (CG) methods are of great importance for solving the large-scale unconstrained optimization problem
\begin{equation}\label{eq:1.1}
\mathop {\min }\limits_{x \in {R^n}} f(x),
\end{equation}
where $f:{R^n} \to R$ is a continuously differentiable function. The key features of CG methods are that they do not require matrix storage. The iterations $\{ {x_n}\} $ satisfy the iterative form
\begin{equation}\label{eq:1.2}
{x_{k + 1}} = {x_k} + {\alpha _k}{d_k}{\rm{,}}
\end{equation} 
where ${\alpha _k}$ is the stepsize and ${d_k}$ is the search direction defined by
\begin{equation}\label{eq:1.3}
{d_{k + 1}} = \left\{ \begin{array}{l}
- {g_{k + 1}},{\rm{                      \quad\quad \quad\; if  }}k = 0,\\
- {g_{k + 1}} + {\beta _k}{d_k},{\rm{            if  }}k > 0,
\end{array} \right.
\end{equation}
where ${g_{k + 1}} = \nabla f({x_{k + 1}})$ and ${\beta _k} \in R$ is called the CG parameter.

For general nonlinear functions, various choices of ${\beta _k}$ cause different CG methods. Some well-known options for ${\beta _k}$ are called FR \cite{22.}, HS \cite{31.}, PRP \cite{40.}, DY \cite{17.} and HZ \cite{28.} formula, and are given by
\begin{align*}
\beta _k^{FR} = \frac{{{{\left\| {{g_{k + 1}}} \right\|}^2}}}{{{{\left\| {{g_k}} \right\|}^2}}},\beta _k^{HS} = \frac{{g_{k + 1}^T{y_k}}}{{d_k^T{y_k}}},\beta _k^{PRP} = \frac{{g_{k + 1}^T{y_k}}}{{{{\left\| {{g_k}} \right\|}^2}}},\beta _k^{DY} = \frac{{{{\left\| {{g_{k + 1}}} \right\|}^2}}}{{d_k^T{y_k}}},
\end{align*}
and
\begin{align*}
\beta _k^{HZ} = \frac{1}{{d_k^T{y_k}}}{\left( {{y_k} - 2{d_k}\frac{{{{\left\| {{y_k}} \right\|}^2}}}{{d_k^T{y_k}}}} \right)^T}{g_{k + 1}},
\end{align*}
where ${y_k} = {g_{k + 1}} - {g_k}$ and $\left\| . \right\|$ denotes the Euclidean norm. Recently, other efficient CG methods have been proposed by different ideas, which can be seen in \cite{16.,21.,28.,29.,32.,41.,42.,50.}. 

With the increasing scale of optimization problems, subspace methods have become a class of very efficient numerical methods because it is not necessary to solve large-scale subproblems at each iteration \cite{80.}. Yuan and Stoer \cite{49.} first put forward the subspace minimization conjugate gradient (SMCG) method, the search direction of which is computed by solving the following problem:
\begin{equation}\label{eq:1.9}
\mathop {\min }\limits_{d \in {\Omega _{k + 1}}} {m_{k + 1}}(d) = g_{k + 1}^Td + \frac{1}{2}{d^T}{B_{k + 1}}d,
\end{equation}
where ${\Omega _{k + 1}} = \left\{ {{g_{k + 1}},{s_k}} \right\}$ and the direction $d$ is given by
\begin{equation}\label{eq:1.10}
d = \mu {g_{k + 1}} + v{s_k},
\end{equation}
where ${B_{k + 1}}$ is an approximation to Hessian matrix, $\mu $ and $v$ are parameters and ${s_k} = {x_{k + 1}} - {x_k}.$ The detailed information of subspace technique can be referred to \cite{1.,30.,33.,46.,51.}. The SMCG method can be considered as a generalization of CG method and it reduces to the linear CG method when it uses the exact line search condition and objective function is convex quadratic function. Based on the analysis of the SMCG method, Dai and Kou \cite{18.} made a theoretical analysis of the BBCG by combining the Barzilai-Borwein (BB) idea \cite{2.} with the SMCG. Liu and Liu \cite{34.} presented an efficient Barzilai-Borwein conjugate gradient method (SMCG\_BB) with the generalized Wolfe line search for unconstrained optimization. Li, Liu and Liu \cite{54.} deliver a subspace minimization conjugate gradient method based on conic model for unconstrained optimization (SMCG\_Conic).

Generally, the iterative methods are often based on a quadratic model because the quadratic model can approximate the objective function well at a small neighborhood of the minimizer. However, when iterative point is far from the minimizer, the quadratic model might not work well if the objective function possesses high non-linearity \cite{43.,48.}. In theory, the successive gradients generated by the conjugate gradient method applied to a quadratic function should be orthogonal. However, for some ill-conditioned problems, orthogonality is quickly lost due to the rounding errors, and the convergence is much slower than expected \cite{30.}. There are many methods to deal with ill-conditioned problems, among which regularization method is one of the effective methods. Recently, $p-$regularized subproblem plays an important role in more regularization approaches \cite{11.,26.,36.} and some $p - $regularization algorithms for unconstrained optimization enjoy a growing interest \cite{3.,6.,12.,11.}. The idea is to incorporate a local quadratic approximation of the objective function with a weighted regularization term $({\sigma _k}/p){\left\| x \right\|^p},p > 2,$ and then globally minimize it at each iteration. Interestingly, Cartis et al. \cite{11.,12.} proved that, under suitable assumptions, $p - $regularization algorithmic scheme is able to achieve superlinear convergence. The most common choice to regularize the quadratic approximation is $p - $regularization with $p = 3,$ which is known as the cubic regularization, since functions of this form are used as local models (to be minimized) in many algorithmic frameworks for unconstrained optimization \cite{4.,5.,6.,9.,10.,11.,12.,20.,24.,26.,36.,38.,45.}. The cubic regularization was first introduced  by Griewank \cite{26.} and was later considered by many authors with global convergence and complexity analysis, see \cite{12.,36.,45.}. 

Recently, how to approximate the $p - $regularized subproblem solution has become a hot research topic. Practical approaches to get an approximate solution are proposed in \cite{6.,25.}, where the solution of the secular equation is typically approximated over specific evolving subspaces using Krylov methods. The main drawback of such approaches is  the large amount of calculation, because they may need to solve multiple linear systems in turn.

In this paper, motivated by \cite{66.} and \cite{47.}, the $p - $regularization with a special scaled norm is analyzed and solutions of the new $p - $regularization that arise in unconstrained optimization are considered. Based on \cite{66.} we propose a method to solve it by using a special scaled norm in the $p - $regularized subproblem. According to the advantages of the new $p - $regularization method with SMCG method, we propose two new subspace minimization conjugate gradient methods. In our algorithms, if the objective function is close to a quadratic, we use a quadratic approximation model in a two-dimensional subspace to generate the direction; otherwise, $p - $regularization model is considered. We prove that the search direction possesses the sufficient descent property and the proposed methods satisfy the global convergence under mild conditions. We present some numerical results, which show that the proposed methods are very promising.

The remainder of this paper is organized as follows. In Section \ref{sec:2}, we will state the form of $p - $regularized subproblem and provide how to solve the $p - $regularization problem based on the special $p - $regularization model. Four choices of search direction by minimizing the approximate models including $p - $regularization and quadratic model on certain subspace are presented in Section \ref{sec:3}. In Section \ref{sec:4}, we describe two algorithms and discuss some important properties of the search direction in detail. In Section \ref{sec:5}, we establish the convergence of the proposed methods under mild conditions. Some performances of the proposed methods are reported in Section \ref{sec:6}. Conclusions and discussions are presented in the last section.
\section {The $p - $regularized Subproblem}
\label{sec:2}
In this section, we will briefly introduce several forms of the $p - $regularized subproblem by using a special scaled norm and provide the solutions of the resulting problems in the whole space and the two-dimensional subspace, respectively. The chosen scaled norm is of the form ${\left\| x \right\|_A} = \sqrt {{x^T}Ax} ,$ where $A$ is a symmetric positive definite matrix. After analysis, we will mainly consider two special cases: (I) $A$ is the Hessian matrix. In this case, the $p-$regularized subproblem has the unique solution; (II) $A$ is the identity matrix. In this case, the $p-$regularized subproblem is the same as the general form.
\subsection{ The Form in the Whole Space}
\label{sec:2.1}
The general form of the $p - $regularized subproblem is: 
\begin{equation}\label{eq:2.1}
\mathop {\min }\limits_{x \in {R^n}} h(x) = {c^T}x + \frac{1}{2}{x^T}Hx + \frac{\sigma }{p}{\left\| x \right\|^p},
\end{equation}
where $p>2,$ $c \in {R^n},$ $\sigma  > 0$ and $H \in {R^{n \times n}}$ is a symmetric matrix.

As for how to solve the above problem, the following theorem is given.\\
\textbf{Theorem 2.1} [ \cite{47.}, Thm.1.1 ] The point ${x^*}$ is a global minimizer of (\ref{eq:2.1}) if and only if
\begin{equation}\label{eq:2.2}
(H + \sigma {\left\| {{x^*}} \right\|^{p - 2}}I){x^*} =  - c, \;\; H + \sigma {\left\| {{x^*}} \right\|^{p - 2}}I  \succeq   0.
\end{equation}
Moreover, the ${l_2}$ norms of all the global minimizers are equal.

Now, we give another form of the $p - $regularized subproblem with a special scaled norm:
\begin{equation}\label{eq:2.3}
\mathop {\min }\limits_{x \in {R^n}} h(x) = {c^T}x + \frac{1}{2}{x^T}Hx + \frac{\sigma }{p}\left\| x \right\|_A^p,
\end{equation}
where $A \in {R^{n \times n}}$ is a symmetric positive definite matrix. 

By setting $y = {A^{\frac{1}{2}}}x,$ (\ref{eq:2.3}) can be arranged as follows:
\begin{equation}\label{eq:2.4}
\mathop {\min }\limits_{y \in {R^n}} h(y) = {({A^{ - \frac{1}{2}}}c)^T}y + \frac{1}{2}{y^T}{A^{ - \frac{1}{2}}}H{A^{ - \frac{1}{2}}}y + \frac{\sigma }{p}{\left\| y \right\|^p}.
\end{equation}
According to Theorem 2.1, we know that the point ${y^*}$ is a global minimizer of (\ref{eq:2.4}) if and only if
\begin{equation}\label{eq:2.5}
({A^{ - \frac{1}{2}}}H{A^{ - \frac{1}{2}}} + \sigma {\left\| {{y^*}} \right\|^{p - 2}}I){y^*} =  - {A^{ - \frac{1}{2}}}c,
\end{equation}
\begin{equation}\label{eq:2.6}
{A^{ - \frac{1}{2}}}H{A^{ - \frac{1}{2}}} + \sigma {\left\| {{y^*}} \right\|^{p - 2}}I \succeq 0.
\end{equation}

Let $V \in {R^{n \times n}}$ be an orthogonal matrix such that
\begin{align*}
{V^T}({A^{ - \frac{1}{2}}}H{A^{ - \frac{1}{2}}})V = Q,
\end{align*}
where $Q = dia{g_{i = 1, \cdots ,n}}\{ {\mu _i}\} $ and ${\mu _1} \le {\mu _2} \le  \cdots  \le {\mu _n}$ are the eigenvalues of ${A^{ - \frac{1}{2}}}H{A^{ - \frac{1}{2}}}.$ Now we can introduce the vector $a \in {R^n}$ such that
\begin{equation}\label{eq:2.7}
y = Va.
\end{equation}

Denote $z = \left\| y \right\|$ and pre-multiplying (\ref{eq:2.5}) by ${V^T},$ we get
\begin{equation}\label{eq:2.8}
(Q + \sigma {z^{p - 2}}I)a =  - \beta ,
\end{equation}
where $\beta  = {V^T}({A^{ - \frac{1}{2}}}c).$

The expression (\ref{eq:2.8}) can be equivalently written as
\begin{align*}
{a_i} = \frac{{ - {\beta _i}}}{{{\mu _i} + \sigma {z^{p - 2}}}}, i = 1,2, \cdots ,n,
\end{align*}
where ${a_i}$ and ${\beta _i}$ are the components of vectors $a$ and $\beta $, respectively. By the way, if ${\mu _i} + \sigma {z^{p - 2}} = 0,$ it means $\beta  = 0$ from (\ref{eq:2.8}).

From (\ref{eq:2.7}), we have an equation about $z$:
\begin{equation}\label{eq:2.9}
{z^2} = {y^T}y = {a^T}a = \mathop \sum \limits_{i = 1}^n \frac{{\beta _i^2}}{{{{({\mu _i} + \sigma {z^{p - 2}})}^2}}}.
\end{equation}

Denote
\begin{align*}
\phi \left( z \right) = \sum\limits_{i = 1}^n {\frac{{\beta _i^2}}{{{{({\mu _i} + \sigma {z^{p - 2}})}^2}}}}  - {z^2}.
\end{align*}
We can easily obtain
\begin{align*}
\phi '\left( z \right) = \sum\limits_{i = 1}^n {\frac{{ - 2\sigma \left( {p - 2} \right)\beta _i^2{z^{p - 3}}\left( {{\mu _i} + \sigma {z^{p - 2}}} \right)}}{{{{\left( {{\mu _i} + \sigma {z^{p - 2}}} \right)}^4}}}}  - 2z.
\end{align*}
It follows from $p > 2,$ $z > 0$ and $\sigma  > 0$ that $\phi '\left( z \right) < 0$, which indicates that $\phi \left( z \right)$ is monotonically decreasing in the interval $\left[ {0, + \infty } \right).$ Moreover, we can observe that $\phi \left( 0 \right) > 0,$ when $\beta  \ne 0,$ and $\mathop {\lim }\limits_{z \to \infty } \phi \left( z \right) =  - \infty .$   So, there exists a unique positive solution to (\ref{eq:2.9}) when $\beta  \ne 0.$ On the other hand, if $\beta  = 0,$ $z = 0$ is the only solution of (\ref{eq:2.9}) in which means ${x^*} = 0$ is the only global minimizer of (\ref{eq:2.3}).

Based on the above derivation and analysis, we can get the following theorem.\\
\textbf{Theorem 2.2} The point ${x^*}$ is a global minimizer of (\ref{eq:2.3}) if and only if
\begin{equation}\label{eq:2.10}
\left( {H + \sigma {{\left( {{z^*}} \right)}^{p - 2}}A} \right){x^*} =  - c,
\end{equation}
\begin{equation}\label{eq:2.11}
H + \sigma {\left( {{z^*}} \right)^{p - 2}}A \succeq 0,
\end{equation}
where ${z^*}$ is the unique non-negative root of the equation
\begin{equation}\label{eq:2.12}
{z^2} = \mathop \sum \limits_{i = 1}^n \frac{{\beta _i^2}}{{{{({\mu _i} + \sigma {z^{p - 2}})}^2}}}.
\end{equation}
Moreover, the ${l_A}$ norms of all the global minimizers are equal.

Now, let us consider a special case that $H \succ 0$ and $A = H.$ It is clear that $H + \sigma {z^{p - 2}}H$ is always a positive definite matrix since $\sigma  > 0$ and $z \ge 0.$ So, the global minimizer of (\ref{eq:2.3}) is unique.\\
\textbf{Inference 2.3} Let $H \succ 0,$ $A = H,$ then the point ${x^*} = \frac{{ - 1}}{{1 + \sigma {{\left( {{z^*}} \right)}^{p - 2}}}}{H^{ - 1}}c$ is the only global minimizer of (\ref{eq:2.3}) and ${z^*}$ is the unique non-negative solution to the equation
\begin{equation}\label{eq:2.13}
\sigma {z^{p - 1}} + z - \sqrt {{c^T}{H^{ - 1}}c}  = 0.
\end{equation}
\textbf{Remark 1} i) $c = 0.$
	It is obvious that the equation (\ref{eq:2.13}) becomes
	\begin{align*}
	\sigma {z^{p - 1}} + z = 0,
	\end{align*}
	that is 
	\begin{align*}
	z\left( {\sigma {z^{p - 2}} + 1} \right) = 0.
	\end{align*}
	From $\sigma  > 0,$ we know ${z^*} = 0$ is the unique non-negative solution to the equation (\ref{eq:2.13}).\\
	ii) $c \ne 0.$
Denote
\begin{equation}\label{eq:2.14}
\psi (z) = \sigma {z^{p - 1}} + z - \sqrt {{c^T}{H^{ - 1}}c}.
\end{equation}
We can easily obtain
\begin{align*}
\psi '(z) = \sigma (p - 1){z^{p - 2}} + 1 > 0,
\end{align*}
which indicates that the $\psi (z)$ is monotonically increasing. From $\psi (0) < 0$ and $\psi \left( {\sqrt {{c^T}{H^{ - 1}}c} } \right) > 0,$ we know that ${z^*}$ is the unique positive solution to the equation (\ref{eq:2.13}). 

\subsection{The Form in the Two-Dimensional Space}
\label{sec:2.2}
Let $g$ and $s$ be two linearly independent vectors. Denote $\Omega  = \left\{ {\left. d \right|d = \mu g + \nu s,\mu ,\nu  \in R} \right\}.$  In this part, we suppose that $H$ is symmetric and positive definite and $y=Hs.$  

We consider the following problem
\begin{equation}\label{eq:2.15}
\mathop {\min }\limits_{d \in \Omega } h(d) = {c^T}d + \frac{1}{2}{d^T}Hd + \frac{\sigma }{p}\left\| d \right\|_A^p.
\end{equation}
Obviously, when $A = H,$ problem (\ref{eq:2.15}) can be translated into
\begin{equation}\label{eq:2.18}
	\mathop {\min }\limits_{\mu ,\nu  \in R} {\left( {\begin{array}{*{20}{c}}
				{{g^T}c}\\
				{{s^T}c}
		\end{array}} \right)^T}\left( {\begin{array}{*{20}{c}}
			\mu \\
			\nu 
	\end{array}} \right) + \frac{1}{2}{\left( {\begin{array}{*{20}{c}}
				\mu \\
				\nu 
		\end{array}} \right)^T}B\left( {\begin{array}{*{20}{c}}
			\mu \\
			\nu 
	\end{array}} \right) + \frac{\sigma }{p}\left\| {\left( {\begin{array}{*{20}{c}}
				\mu \\
				\nu 
		\end{array}} \right)} \right\|_B^p,
\end{equation}
where $\rho  = {g^T}Hg$, and $B = \left( {\begin{array}{*{20}{c}}
	\rho &{{g^T}y}\\
	{{g^T}y}&{{y^T}s}
	\end{array}} \right)$ is a symmetric and positive definite matrix since the $H$ is a symmetric positive definite matrix and the two vectors $g$ and $s$ are linearly independent.

By the Inference 2.3, we can obtain the unique solution of (\ref{eq:2.18}):
\begin{equation}\label{eq:2.19}
\left( {\begin{array}{*{20}{c}}
	{{\mu ^*}}\\
	{{\nu ^*}}
	\end{array}} \right) = \frac{{ - 1}}{{1 + \sigma {{\left( {{z^*}} \right)}^{p - 2}}}}{B^{ - 1}}\left( {\begin{array}{*{20}{c}}
	{{g^T}c}\\
	{{s^T}c}
	\end{array}} \right),
\end{equation}
where ${z^*}$ is the unique non-negative solution to $\sigma {z^{p - 1}} + z - \sqrt {{{\left( {\begin{array}{*{20}{c}}
					{{g^T}c}\\
					{{s^T}c}
					\end{array}} \right)}^T}{B^{ - 1}}\left( {\begin{array}{*{20}{c}}
			{{g^T}c}\\
			{{s^T}c}
			\end{array}} \right)}  = 0.$

When $A = I,$ we obtain from (\ref{eq:2.15}) that
\begin{equation}\label{eq:2.20}
\mathop {\min }\limits_{\mu ,\nu  \in R} {\left( {\begin{array}{*{20}{c}}
		{{g^T}c}\\
		{{s^T}c}
		\end{array}} \right)^T}\left( {\begin{array}{*{20}{c}}
	\mu \\
	\nu 
	\end{array}} \right) + \frac{1}{2}{\left( {\begin{array}{*{20}{c}}
		\mu \\
		\nu 
		\end{array}} \right)^T}B\left( {\begin{array}{*{20}{c}}
	\mu \\
	\nu 
	\end{array}} \right) + \frac{\sigma }{p}\left\| {\left( {\begin{array}{*{20}{c}}
		\mu \\
		\nu 
		\end{array}} \right)} \right\|_E^p,
\end{equation}
where $E = \left( {\begin{array}{*{20}{c}}
	{{{\left\| g \right\|}^2}}&{{g^T}s}\\
	{{g^T}s}&{{{\left\| s \right\|}^2}}
	\end{array}} \right)$ is  positive definite due to the linear independence of vectors $g$ and $s.$

By the Theorem 2.2, we can gain the unique solution to (\ref{eq:2.20}):
\begin{equation}\label{eq:2.21}
\left( {\begin{array}{*{20}{c}}
	{{\mu ^*}}\\
	{{\nu ^*}}
	\end{array}} \right) =  - {\left( {B + \sigma {{\left( {{z^*}} \right)}^{p - 2}}E} \right)^{ - 1}}\left( {\begin{array}{*{20}{c}}
	{{g^T}c}\\
	{{s^T}c}
	\end{array}} \right),
\end{equation}
where ${z^*}$ is the unique non-negative solution to (\ref{eq:2.12}) in which $0 < {\mu _1} \le {\mu _2}$ are the eigenvalues of ${E^{ - \frac{1}{2}}}B{E^{ - \frac{1}{2}}},$ $\beta  = {V^T}\left( {{E^{ - \frac{1}{2}}}\left( {\begin{array}{*{20}{c}}
		{{g^T}c}\\
		{{s^T}c}
		\end{array}} \right)} \right).$
\section {The Search Direction and The Initial Stepsize}
\label{sec:3}
In this section, based on the different choices of special scaled norm, we derive two new directions by minimizing the two $p - $regularization models of the objective function on the subspace ${\Omega _k} = span\left\{ {{g_k},{s_{k - 1}}} \right\}.$ The selection criteria for how to choose the initial stepsize is given. For the rest, we assume that $s_k^T{y_k} > 0$ guaranteed by the condition (\ref{eq:1.5}).
\subsection{Derivation of The New Search Direction}
\label{sec:3.1}
The parameter ${t_k}$ by Yuan \cite{50.} is used describe how $f(x)$ is close to a quadratic function on the line segment between ${x_{k - 1}}$ and ${x_k},$ and defined by
\begin{equation}\label{eq:3.1}
{t_k} = \left| {\frac{{2\left( {{f_{k - 1}} - {f_k} + g_k^T{s_{k - 1}}} \right)}}{{s_{k - 1}^T{y_{k - 1}}}} - 1} \right|.
\end{equation}

On the other hand, the ratio
\begin{equation}\label{eq:3.3}
{\theta _k} = \frac{{{f_{k - 1}} - {f_k}}}{{0.5s_{k - 1}^T{y_{k - 1}} - g_k^T{s_{k - 1}}}}
\end{equation}
shows difference between the actual reduction and the predicted reduction for the quadratic model.

If the following condition \cite{39.} holds, namely,
\begin{equation}\label{eq:3.2}
{t_k} \le {c_1}\;\;{\rm{or}}\;\;\left( {{t_k} \le {c_2}\;\;{\rm{and}}\;\;{t_{k - 1}} \le {c_2}} \right)
\end{equation}
or
\begin{equation}\label{eq:3.4}
\left| {{\theta _k} - 1} \right| < \gamma,  
\end{equation}
where ${c_1},{c_2}{\;\;\rm{ and \;\;}}\gamma $ are small positive constants, then $f(x)$ might be very close to a quadratic on the line segment between ${x_{k - 1}}$ and ${x_k}.$ We choose the quadratic model. 

Moreover, if the conditions \cite{54.}
\begin{equation}\label{eq:4.1}
{(s_k^T{y_k})^2} \le {10^{ - 5}}{\left\| {{s_k}} \right\|^2}{\left\| {{y_k}} \right\|^2}\;\;{\rm{and}}\;\;{({f_{k + 1}} - {f_k} - 0.5(g_k^T{s_k} + g_{k + 1}^T{s_k}))^2} \le {10^{ - 6}}{\left\| {{s_k}} \right\|^2}{\left\| {{y_k}} \right\|^2}
\end{equation}
hold, then the problem might have very large condition number, which seems to be ill-conditioned. And the current iterative point is far away from the minimizer of problem. At this point, the information might be inaccurate, then we also choose the quadratic model to derive a search direction.

General iterative methods, which are often based on a quadratic model, have been quite successful for solving unconstrained optimization problems, since the quadratic model can approximate the objective function $f\left( x \right)$ well at a small neighborhood of ${x_k}$ in many cases. Consequently, when the condition (\ref{eq:3.2}), (\ref{eq:3.4}) or (\ref{eq:4.1})  holds, the quadratic approximation model (\ref{eq:1.9}) is preferable. However, when the conditions (\ref{eq:3.2}), (\ref{eq:3.4}) and (\ref{eq:4.1}) do not hold, the iterative point is far away from the minimizer, the quadratic model may not very well approximate the original problem. Thus in this case,  we select the $p - $regularization model which could include more useful information of the objective function to approximate the original problem.

For general functions, if the condition
	\begin{equation}\label{eq:3.18}
	{\xi _1} \le \frac{{s_{k - 1}^T{y_{k - 1}}}}{{{{\left\| {{s_{k - 1}}} \right\|}^2}}} \le \frac{{{{\left\| {{y_{k - 1}}} \right\|}^2}}}{{s_{k - 1}^T{y_{k - 1}}}} \le {\xi _2}
	\end{equation}
	holds, where ${\xi _1}$ and ${\xi _2}$ are positive constants, then the condition number of the Hessian matrix might be not very large. In this case, we consider the quadratic approximation model or the $p-$regularization model. 

Now we divide it into following four cases to derive the search direction.\\
\textbf{Case 1.} When the condition (\ref{eq:3.18}) holds and any of the conditions (\ref{eq:3.2}, \ref{eq:3.4}, \ref{eq:4.1}) do not hold, we consider the following $p - $regularized subproblem
\begin{equation}\label{eq:3.5}
\mathop {\min }\limits_{{d_k} \in {\Omega _k}} {m_k}\left( {{d_k}} \right) = d_k^T{g_k} + \frac{1}{2}d_k^T{{H}_k}{d_k} + \frac{1}{p}{\sigma _k}\left\| {{d_k}} \right\|_{{{A}_k}}^p,
\end{equation}
where ${{H}_k}$ is a symmetric and positive definite approximation to Hessian matrix satisfying the equation ${{H}_k}{s_{k - 1}} = {y_{k - 1}},$ ${{A}_k}$ is a symmetric positive definite matrix, ${\sigma _k}$ is a dynamic non-negative regularization parameter and ${\Omega _k} = span\{ {g_k},{s_{k - 1}}\} .$

Denote
\begin{equation}\label{eq:3.6}
{d_k} = {\mu _k}{g_k} + {\nu _k}{s_{k - 1}},
\end{equation}
where ${\mu _k}$ and ${\nu _k}$ are parameters to be determined.

In the following, we will discuss that ${{A}_k} = {{H}_k}$ and ${{A}_k} = I$ in two parts.

(I) ${{A}_k} = {{H}_k}$

It is easy to see the problem (\ref{eq:3.5}) is similar to the problem (\ref{eq:2.18}), we obtain
\begin{equation}\label{eq:3.7}
\mathop {\min }\limits_{{\mu _k},{\nu _k} \in R} {\left( {\begin{array}{*{20}{c}}
		{{{\left\| {{g_k}} \right\|}^2}}\\
		{g_k^T{s_{k - 1}}}
		\end{array}} \right)^T}\left( {\begin{array}{*{20}{c}}
	{{\mu _k}}\\
	{{\nu _k}}
	\end{array}} \right) + \frac{1}{2}{\left( {\begin{array}{*{20}{c}}
		{{\mu _k}}\\
		{{\nu _k}}
		\end{array}} \right)^T}{B_k}\left( {\begin{array}{*{20}{c}}
	{{\mu _k}}\\
	{{\nu _k}}
	\end{array}} \right) + \frac{{{\sigma _k}}}{p}\left\| {\left( {\begin{array}{*{20}{c}}
		{{\mu _k}}\\
		{{\nu _k}}
		\end{array}} \right)} \right\|_{{B_k}}^p,
\end{equation}
where ${\rho _k} \approx g_k^T{{H}_k}{g_k}$ and ${B_k} =  \left( {\begin{array}{*{20}{c}}
	{{\rho _k}}&{g_k^T{y_{k - 1}}}\\
	{g_k^T{y_{k - 1}}}&{s_{k - 1}^T{y_{k - 1}}}
	\end{array}} \right)$.

It is very important for how to choose the two parameters ${{{ \rho }_k}}$ and ${{\sigma _k}}$ in (\ref{eq:3.7}).

Motivated by the Barzilai-Borwein method, Dai and Kou \cite{18.} proposed a BBCG3 method with the very efficient parameter $\rho _k^{BBCG3} = \frac{3}{2}\frac{{{{\left\| {{y_{k - 1}}} \right\|}^2}}}{{s_{k - 1}^T{y_{k - 1}}}}{\left\| {{g_k}} \right\|^2}$ and considered it a good estimation of the $g_k^T{{H}_k}{g_k}$. So in this paper, we choose ${\rho _k} = \rho _k^{BBCG3}$ in the above function that will make  ${B_k}$ positive, which guarantees definite the unique solution to (\ref{eq:3.7}).

There are many ways \cite{11.,24.} to get the value of ${\sigma _k}$, and the interpolation condition is one of them. Here, we use interpolation condition to get it. By imposing the following interpolation condition: 
\begin{align*}
{f_{k - 1}} = {f_k} - g_k^T{s_{k - 1}} + \frac{1}{2}s_{k - 1}^T{y_{k - 1}} + \frac{{{\sigma _k}}}{p}{(s_{k - 1}^T{y_{k - 1}})^{\frac{p}{2}}},
\end{align*}
we obtain 
\begin{align*}
{\sigma _k} = \frac{{p({f_{k - 1}} - {f_k} + g_k^T{s_{k - 1}} - \frac{1}{2}s_{k - 1}^T{y_{k - 1}})}}{{{{(s_{k - 1}^T{y_{k - 1}})}^{\frac{p}{2}}}}}.
\end{align*}
In order to ensure that ${\sigma _k} \ge 0,$ we set
\begin{align*}
 {\sigma _k} = \frac{{p\left| {{f_{k - 1}} - {f_k} + g_k^T{s_{k - 1}} - \frac{1}{2}s_{k - 1}^T{y_{k - 1}}} \right|}}{{{{\left( {s_{k - 1}^T{y_{k - 1}}} \right)}^{\frac{p}{2}}}}}.
\end{align*}

From (\ref{eq:2.19}), we can get the unique solution to (\ref{eq:3.7}):
 \begin{equation}\label{eq:3.8}
{\mu _k} = \frac{1}{{\left( {1 + {\sigma _k}{{\left( {{z^*}} \right)}^{p - 2}}} \right){\Delta _k}}}\left( {g_k^T{y_{k - 1}}g_k^T{s_{k - 1}} - s_{k - 1}^T{y_{k - 1}}{{\left\| {{g_k}} \right\|}^2}} \right),
 \end{equation} 
\begin{equation}\label{eq:3.9}
{\nu _k} = \frac{1}{{\left( {1 + {\sigma _k}{{\left( {{z^*}} \right)}^{p - 2}}} \right){\Delta _k}}}\left( {g_k^T{y_{k - 1}}{{\left\| {{g_k}} \right\|}^2} - {\rho _k}g_k^T{s_{k - 1}}} \right),
\end{equation}
where ${\Delta _k} = \left| {\begin{array}{*{20}{c}}
	{{\rho _k}}&{g_k^T{y_{k - 1}}}\\
	{g_k^T{y_{k - 1}}}&{s_{k - 1}^T{y_{k - 1}}}
	\end{array}} \right| = {\rho _k}s_{k - 1}^T{y_{k - 1}} - {(g_k^T{y_{k - 1}})^2} > 0$ and ${z^*}$ is the unique positive solution to 
\begin{equation}\label{eq:3.10}
{\sigma _k}{z^{p - 1}} + z - \sqrt {{{\left( {\begin{array}{*{20}{c}}
				{{{\left\| {{g_k}} \right\|}^2}}\\
				{g_k^T{s_{k - 1}}}
				\end{array}} \right)}^T}{B_k}^{ - 1}\left( {\begin{array}{*{20}{c}}
		{{{\left\| {{g_k}} \right\|}^2}}\\
		{g_k^T{s_{k - 1}}}
		\end{array}} \right)}  = 0.
\end{equation}

We denote $\tilde q = \sqrt {\left( {\begin{array}{*{20}{c}}
		{{{\left\| {{g_k}} \right\|}^2}}\\
		{g_k^T{s_{k - 1}}}
		\end{array}} \right)B_k^{ - 1}\left( {\begin{array}{*{20}{c}}
		{{{\left\| {{g_k}} \right\|}^2}}\\
		{g_k^T{s_{k - 1}}}
		\end{array}} \right)}.$ Substituting $\tilde q$ into (\ref{eq:3.10}), we get
\begin{equation}\label{eq:3.11}
{\sigma _k}{z^{p - 1}} + z - \tilde q = 0.
\end{equation}
Since it is difficult to obtain the exact root of (\ref{eq:3.11}) when $p$ is large, we only consider $p=3$ and $p=4$ for simplicity. 

(i) $p = 3.$ It is not difficult to know the unique positive solution to (\ref{eq:3.11})
\begin{equation}\label{eq:3.12}
{z^*} = \frac{{2\tilde q}}{{1 + \sqrt {1 + 4{\sigma _k}\tilde q} }}.
\end{equation}

(ii) $p = 4.$ According to the formula of extracting roots on cubic equation and $z > 0,$ the unique positive solution to (\ref{eq:3.11}) can be obtained
\begin{equation}\label{eq:3.13}
{z^*} = \sqrt[3]{{\frac{{\tilde q}}{{2{\sigma _k}}} + \sqrt {\frac{{{{\tilde q}^2}}}{{4\sigma _k^2}} + {{\left( {\frac{1}{{3{\sigma _k}}}} \right)}^3}} }} + \sqrt[3]{{\frac{{\tilde q}}{{2{\sigma _k}}} - \sqrt {\frac{{{{\tilde q}^2}}}{{4\sigma _k^2}} + {{\left( {\frac{1}{{3{\sigma _k}}}} \right)}^3}} }}.
\end{equation}

For ensuring the sufficient descent condition of the direction produced by (\ref{eq:3.8}) and (\ref{eq:3.9}), if ${\sigma _k}{\left( {{z^*}} \right)^{p - 2}} > 1,$ we set ${\sigma _k}{\left( {{z^*}} \right)^{p - 2}} = 1,$ where ${{z^*}}$ is determined by (\ref{eq:3.12}) or (\ref{eq:3.13}).

(II) ${{A}_k} = I$ 

Based on the analysis of (I), we can get the following problem similarly:
\begin{equation}\label{eq:3.14}
\mathop {\min }\limits_{{\mu _k},{\nu _k} \in R} {\left( {\begin{array}{*{20}{c}}
		{{{\left\| {{g_k}} \right\|}^2}}\\
		{g_k^T{s_{k - 1}}}
		\end{array}} \right)^T}\left( {\begin{array}{*{20}{c}}
	{{\mu _k}}\\
	{{\nu _k}}
	\end{array}} \right) + \frac{1}{2}{\left( {\begin{array}{*{20}{c}}
		{{\mu _k}}\\
		{{\nu _k}}
		\end{array}} \right)^T}{B_k}\left( {\begin{array}{*{20}{c}}
	{{\mu _k}}\\
	{{\nu _k}}
	\end{array}} \right) + \frac{{{\sigma _k}}}{p}\left\| {\left( {\begin{array}{*{20}{c}}
		{{\mu _k}}\\
		{{\nu _k}}
		\end{array}} \right)} \right\|_{{E_k}}^p,
\end{equation}
where ${E_k} = \left( {\begin{array}{*{20}{c}}
	{{{\left\| {{g_k}} \right\|}^2}}&{g_k^T{s_{k - 1}}}\\
	{g_k^T{s_{k - 1}}}&{{{\left\| {{s_{k - 1}}} \right\|}^2}}
	\end{array}} \right)$ and ${\rho _k}$, ${B_k}$ are the same as those in problem (\ref{eq:3.7}).

Similarly, we still use the interpolation condition to determine ${\sigma _k}$:  
\begin{align*}
{f_{k - 1}} = {f_k} - g_k^T{s_{k - 1}} + \frac{1}{2}s_{k - 1}^T{y_{k - 1}} + \frac{{{\sigma _k}}}{p}{\left\| {{s_{k - 1}}} \right\|^{\frac{p}{2}}},
\end{align*}
we get
\begin{align*}
	{\sigma _k} = \frac{{p\left| {{f_{k - 1}} - {f_k} + g_k^T{s_{k - 1}} - \frac{1}{2}s_{k - 1}^T{y_{k - 1}}} \right|}}{{{{\left\| {{s_{k - 1}}} \right\|}^{\frac{p}{2}}}}}.
	\end{align*}
 
According to (\ref{eq:2.21}), the unique solution to (\ref{eq:3.14}) can be obtained:
\begin{equation}\label{eq:3.15}
{{\hat \mu }_k} = \frac{1}{{{{\bar \Delta }_k}}}\left( {g_k^T{y_{k - 1}}g_k^T{s_{k - 1}} - s_{k - 1}^T{y_{k - 1}}{{\left\| {{g_k}} \right\|}^2} + \lambda {{\left( {g_k^T{s_{k - 1}}} \right)}^2} - \lambda {{\left\| {{s_{k - 1}}} \right\|}^2}{{\left\| {{g_k}} \right\|}^2}} \right),
\end{equation}
\begin{equation}\label{eq:3.16}
{{\hat \nu }_k} = \frac{1}{{{{\bar \Delta }_k}}}\left( {g_k^T{y_{k - 1}}{{\left\| {{g_k}} \right\|}^2} - {\rho _k}g_k^T{s_{k - 1}}} \right),
\end{equation}
where 
\begin{equation}\label{eq:3.17}
{{\bar \Delta }_k} = ({\rho _k} + \lambda {\left\| {{g_k}} \right\|^2})(s_{k - 1}^T{y_{k - 1}} + \lambda {\left\| {{s_{k - 1}}} \right\|^2}) - {(g_k^T{y_{k - 1}} + \lambda g_k^T{s_{k - 1}})^2},
\end{equation}
\begin{align*}
\lambda  = {\sigma _k}{\left( {{z^*}} \right)^{p - 2}}
\end{align*}
and ${z^*} $ satisfies the equation (\ref{eq:2.12}), which can be solved by tangent method \cite{37.}. For ensuring the sufficient descent of the direction produced by (\ref{eq:3.15}) and (\ref{eq:3.16}), if ${\sigma _k}{\left( {{z^*}} \right)^{p - 2}} > \frac{{{{\left\| {{y_{k - 1}}} \right\|}^2}}}{{s_{k - 1}^T{y_{k - 1}}}},$ we set $\lambda  = \frac{{{{\left\| {{y_{k - 1}}} \right\|}^2}}}{{s_{k - 1}^T{y_{k - 1}}}}.$\\
\textbf{Remark 2} It is worth emphasizing that in the process of finding the direction, $\left( {\begin{array}{*{20}{c}}
		{{{\left\| {{g_k}} \right\|}^2}}\\
		{g_k^T{s_{k - 1}}}
		\end{array}} \right) \ne 0$, which is equivalent to the problem (\ref{eq:2.3}) in which $c \ne 0$.\\
\textbf{Case 2.} When the condition (\ref{eq:3.18}) holds and one of the conditions (\ref{eq:3.2}, \ref{eq:3.4}, \ref{eq:4.1}) at least holds, we choose the quadratic model which corresponds to (\ref{eq:3.7}) with $\sigma_k=0.$ So the parameters in (\ref{eq:3.6}) are generated by solving  (\ref{eq:3.8}) and (\ref{eq:3.9})  with ${\sigma _k} = 0:$ \\ 
\begin{equation}\label{eq:3.21}
{{\bar \mu }_k} = \frac{1}{{{\Delta _k}}}(g_k^T{y_{k - 1}}g_k^T{s_{k - 1}} - s_{k - 1}^T{y_{k - 1}}{\left\| {{g_k}} \right\|^2}),
\end{equation}
\begin{equation}\label{eq:3.22}
{{\bar \nu }_k} = \frac{1}{{{\Delta _k}}}(g_k^T{y_{k - 1}}{\left\| {{g_k}} \right\|^2} - {\rho _k}g_k^T{s_{k - 1}}).
\end{equation} 
\textbf{Case 3.} If the exact line search is used, the direction in Case 2 is parallel to the HS direction with convex quadratic functions. It is known that the conjugate condition, namely, $d_{k + 1}^T{y_k} = 0,$ still holds whether the line search is exact or not for HS conjugate gradient method. 

If the condition (\ref{eq:3.18}) does not hold and the conditions
\begin{equation}\label{eq:3.23}
\frac{{\left| {g_k^T{y_{k - 1}}g_k^T{s_{k - 1}}} \right|}}{{s_{k - 1}^T{y_{k - 1}}{{\left\| {{g_k}} \right\|}^2}}} \le {\xi _3}{\rm{  \;\;and\;\;  }}{\xi _1} \le \frac{{s_{k - 1}^T{y_{k - 1}}}}{{{{\left\| {{s_{k - 1}}} \right\|}^2}}}
\end{equation}
hold, where $0 \le {\xi _3} \le 1,$ then ${{\bar \mu }_k}$ in Case 2 is close to -1, then we use the HS conjugate gradient direction. Besides, with the finite-termination property of the HS method for exact convex quadratic programming, such choice of the direction might lead to a rapid convergence rate of our algorithm.\\
\textbf{Case 4.} If the condition  (\ref{eq:3.18}) does not hold and the condition (\ref{eq:3.23}) does not hold, then we choose the negative gradient as the search direction, namely,
\begin{equation}\label{eq:3.24}
{d_k} =  - {g_k}.
\end{equation} 

In conclusion, the new search direction can be stated as
\begin{align*}
{d_k} = \left\{ {\begin{array}{*{20}{c}}
	{{\mu _k}{g_k} + {\nu _k}{s_{k - 1}},}&{}&{{\rm{if}}}&{\left( {\ref{eq:3.18}} \right)}&{{\rm{holds}}}&{{\rm{and}}}&{{\rm{any}}}&{{\rm{of}}}&{\left( {{\ref{eq:3.2}, \ref{eq:3.4}, \ref{eq:4.1}}} \right)}&{{\rm{do}}}&{{\rm{not}}}&{{\rm{hold,}}}\\
	{{{\bar \mu }_k}{g_k} + {{\bar \nu }_k}{s_{k - 1}},}&{}&{{\rm{if}}}&{\left( {\ref{eq:3.18}} \right)}&{{\rm{holds}}}&{{\rm{and}}}&{{\rm{one}}}&{{\rm{of}}}&{\left( {{\ref{eq:3.2}, \ref{eq:3.4}, \ref{eq:4.1}}} \right)}&{{\rm{at}}}&{{\rm{least}}}&{{\rm{holds,}}}\\
	{ - {g_k} + \beta _k^{HS}{d_{k - 1}},}&{}&{{\rm{if}}}&{\left( {\ref{eq:3.18}} \right)}&{{\rm{does}}}&{{\rm{not}}}&{{\rm{hold}}}&{{\rm{and}}}&{\left( {{\ref{eq:3.23}}} \right)}&{{\rm{holds,}}}&{}&{}\\
	{ - {g_k},}&{}&{{\rm{if}}}&{\left( {\ref{eq:3.18}} \right)}&{{\rm{does}}}&{{\rm{not}}}&{{\rm{hold}}}&{{\rm{and}}}&{\left( {{\ref{eq:3.23}}} \right)}&{{\rm{does}}}&{{\rm{not}}}&{{\rm{hold}},}
	\end{array}} \right.
\end{align*}
	where ${\mu _k},{\nu _k}$ are given by (\ref{eq:3.8}), (\ref{eq:3.9}) or (\ref{eq:3.15}), (\ref{eq:3.16}) and ${{\bar \mu }_k}$, ${{\bar \nu }_k}$ are given by (\ref{eq:3.21}), (\ref{eq:3.22}), respectively.\\

\subsection{Choices of The Initial Stepsize and The Wolfe Line Search}
\label{sec:3.2}
It is universally acknowledged that the choice of the initial stepsize and the Wolfe line search are of great importance for an optimization method. In this section, we introduce a strategy to choose the initial stepsize and develop a modified nonmonotone Wolfe line search. 
\subsubsection{Choices of The Initial Stepsize}
Denote 
\begin{align*}
{\phi _k}(\alpha ) = f({x_k} + \alpha {d_k}),\alpha  \ge 0.
\end{align*}
(i) The initial stepsize for the search directions in Case1.-Case3. in Section 3.1.

Similar to \cite{34.}, we choose the initial stepsize as
\begin{align*}
\alpha _k^0 = \left\{ \begin{array}{l}
{{\hat \alpha }_k},\;\;\;\;\;{\rm{    if\;\; (\ref{eq:3.2})\;\; holds\;\; and\;\; }}{{\bar \alpha }_k} > 0,\\
1,\;\;\;\;\;\;\;{\rm{       otherwise,}}
\end{array} \right.
\end{align*} 
where
\begin{align*}
{{\bar \alpha }_k} = \min q({\phi _k}(0),{\phi _k}^\prime (0),{\phi _k}(1)),\;\;{{\hat \alpha }_k} = \min \{ \max\{ {{\bar \alpha }_k},{\lambda _{\min }}\} ,{\lambda _{\max }}\} \;\;{\rm{and}}\;\;{\lambda _{\max }} > {\lambda _{\min }} > 0.
\end{align*}
In the above formula, $q\left( {{\phi _k}\left( 0 \right),{{\phi '}_k}\left( 0 \right),{\phi _k}\left( 1 \right)} \right)$ denotes the interpolation function for the three values ${{\phi _k}\left( 0 \right),}$ ${{{\phi '}_k}\left( 0 \right),}$ and ${{\phi _k}\left( 1 \right).}$ And ${\lambda _{\max }}$ and ${\lambda _{\min }}$ represent two positive parameters.\\
(ii) The initial stepsize for the negative gradient direction (\ref{eq:3.24}).

As we all know, the gradient method with the adaptive BB stepsize \cite{53.} is very efficient for strictly convex quadratic minimization, especially when the condition number is large. In this paper we choose the strategy in \cite{34.}:
\begin{align*}
\alpha _k^0 = \left\{ \begin{array}{l}
\min \{ \max\{ {\tilde{\tilde \alpha}  }_k,{\lambda _{\min }}\} ,{\lambda _{\max }}\} ,{\rm{\;\;if\;\;(\ref{eq:3.2}) \;\;holds,\;\;}}{d_{k - 1}} \ne  - {g_{k - 1}},\;\;{\left\| {{g_k}} \right\|^2} \le 1{\rm{ \;\;and\;\; }}{\tilde{\tilde \alpha}  }_k > 0,\\
{{\bar {\bar \alpha} }_k},{\rm{ \quad \qquad\qquad\qquad\qquad\qquad otherwise}}{\rm{,}}
\end{array} \right.
\end{align*}
where
\begin{align*}
{\bar {\bar \alpha} _k} = \left\{ \begin{array}{l}
\{ \min \{ {\lambda _k}\alpha _k^{B{B_2}},{\lambda _{\max }}\} ,{\lambda _{\min }}\} ,\;\;{\rm{        if}}\;\;g_k^T{s_{k - 1}} > 0,\\
\{ \min \{ {\lambda _k}\alpha _k^{B{B_1}},{\lambda _{\max }}\} ,{\lambda _{\min }}\} ,\;\;{\rm{        if}}\;\;g_k^T{s_{k - 1}} \le 0,
\end{array} \right.,\;\;{{\tilde{\tilde \alpha}  }_k} = \min q({\phi _k}(0),{\phi _k}^\prime (0),{\phi _k}({\bar {\bar \alpha} _k})),
\end{align*}
${\lambda _k}$ is a scaling parameter given by ${\lambda _k} = \left\{ \begin{array}{l}
0.999,\;\;{\rm{    if }}\;\;n > 10{\rm{ \;\;and \;\;Numgra  >  12,}}\\
1,{\rm{\;\;\;\;\;\;\;\;\;otherwise,}}
\end{array} \right.$ \\
where Numgra denotes the number of the successive use of the negative gradient direction.
\subsubsection{Choice of The Wolfe Line Search} 
The line search is an important factor for the overall efficiency of most optimization algorithms. In this paper, we pay attention to the nonmonotone line search proposed by Zhang and Hager \cite{52.} (ZH line search)
\begin{equation}\label{eq:1.4}
	f({x_k} + {\alpha _k}{d_k}) \le {C_k} + \delta {\alpha _k}\nabla f{({x_k})^T}{d_k},
\end{equation}
\begin{equation}\label{eq:1.5}
	\nabla f{({x_k} + {\alpha _k}{d_k})^T}{d_k} \ge \sigma \nabla f{({x_k})^T}{d_k},
\end{equation}
where $0 < \delta  < \sigma  < 1,$ ${C_0} = {f_0},$ ${Q_0} = 1,$ and ${C_k}$ and ${Q_k}$ are updated by
\begin{equation}\label{eq:1.6}
	{Q_{k + 1}} = {\eta _k}{Q_k} + 1,{C_{k + 1}} = \frac{{{\eta _k}{Q_k}{C_k} + f({x_{k + 1}})}}{{{Q_{k + 1}}}},
\end{equation}
where ${\eta _k} \in [0,1].$

It is worth mentioning that some improvements have been made to ZH line search to find a more suitable stepsize and obtain a better convergence result. Specially,
\begin{equation}\label{eq:1.7}
	{C_1} = \min \{ {C_0},{f_1} + 1.0\} , {Q_1} = 2.0,
\end{equation}
when $k \ge 1,$ ${C_{k + 1}}$ and ${Q_{k + 1}}$ are updated by  (\ref{eq:1.6}), where ${\eta _k}$ is taken as
\begin{equation}\label{eq:1.8}
	{\eta _k} = \left\{ \begin{array}{l}
	\eta ,\;\;{\rm{ if }}\bmod (k,l) = 0,\\
	1,{\rm{  \;\;if }}\bmod (k,l) \ne 0,
	\end{array} \right.
\end{equation}
where $l = \max (20,n),$ $\bmod (k,l)$ denotes the remainder for $k$ modulo $l$ and $\eta  = 0.7$ when ${C_k} - {f_{k+1}} > 0.999\left| {{C_k}} \right|,$ otherwise $\eta  = 0.999.$ Such choice of ${\eta _k}$ can be used to control nonmonotonicity dynamically, referred to \cite{35.}.
\section{Algorithms}
\label{sec:4}
In this section, according to the different choices of special scaled norm, we will introduce two new subspace minimization conjugate gradient algorithms based on the $p - $regularization and analyze some theoretical properties of the direction ${d_k}.$

Denote
\begin{align*}
{r_{k - 1}} = \left| {\frac{{{f_k}}}{{{f_{k - 1}} + 0.5(g_{k - 1}^T{s_{k - 1}} + g_k^T{s_{k - 1}})}} - 1} \right|,\;\;{{\bar r}_{k - 1}} = \left| {{f_k} - {f_{k - 1}} - 0.5(g_{k - 1}^T{s_{k - 1}} + g_k^T{s_{k - 1}})} \right|.
\end{align*}
 If ${r_{k - 1}}$ or ${{\bar r}_{k - 1}}$ is close to 0, then the function might be close to a quadratic function. If there are continuously many iterations such that ${r_{k - 1}} \le {\xi _4}$ or ${{\bar r}_{k - 1}} \le {\xi _5},$ where ${\xi _4},{\xi _5} > 0,$ we restart the method with $ - {g_k}$. In addition, if the number of the successive use of CG direction reaches to the threshold MaxRestart, we also restart the method with $ - {g_k}$.

Firstly, we describe the subspace minimization conjugate gradient method in which the direction of the regularization model is generated by the problem (\ref{eq:3.7}), which is called SMCG\_PR1.

\begin{algorithm}[H]
	\caption{ SMCG method with $p - $regularization (SMCG\_PR1)}\label{alg:1}
	\textbf{Step 0.} Given ${x_0} \in {R^n},\;\;\varepsilon  > 0,\;\;0 < \delta  < \sigma  < 1,\;\;{\xi _1},\;\;{\xi _2},\;\;{\xi _3},\;\;{\xi _4},\;\;{\xi _5},\;\;{c_1},\;\;{c_2},\;\;\gamma  \in (0,1),\;\;\alpha _0^{(0)}.$ Let ${C_0} = {f_0},$ ${Q_0} = 1,$ \hspace*{1.1cm}${d_0} =  - {g_0}$ and $k: = 0.$ Set IterRestart :=0, Numgrad :=0, IterQuad :=0, Isnotgra=0,  
	MaxRestart, MinQuad.\\
	\textbf{Step 1.} If ${\left\| {{g_k}} \right\|_\infty } \le \varepsilon ,$ then stop.\\
	\textbf{Step 2.} Compute a stepsize ${\alpha _k} > 0$ satisfying (\ref{eq:1.4}) and (\ref{eq:1.5}). Let ${x_{k + 1}} = {x_k} + {\alpha _k}{d_k}.$ If ${\left\| {{g_k}} \right\|_\infty } \le \varepsilon ,$ then stop. Otherwise, \hspace*{1.1cm} set IterRestart:=IterRestart+1. If ${r_{k - 1}} \le {\xi _4}$ or ${{\bar r}_{k - 1}} \le {\xi _5},$ then IterQuad :=IterQuad+1, else IterQuad :=0.\\
	\textbf{Step 3.} (Calculation of the direction)\\
\hspace*{0.4cm} \textbf{3.1.} If Isnotgra=MaxRestart or (IterQuad=MinQuad  
and IterRestart $ \ne $ IterQuad), then set ${d_{k + 1}} =  - {g_{k + 1}}.$ Set \hspace*{1.1cm} Numgrad := Numgrad+1, Isnotgra :=0 and IterRestart :=0, and go to Step 4. If the condition (\ref{eq:3.18}) holds, go to \hspace*{1.1cm} 3.2;  otherwise go to 3.3.\\
\hspace*{0.4cm} \textbf{3.2.} If the condition (\ref{eq:3.2}) or (\ref{eq:3.4}) or (\ref{eq:4.1}) holds, compute the search direction ${d_{k + 1}}$ by (\ref{eq:3.6}) with (\ref{eq:3.21}) and (\ref{eq:3.22}). \hspace*{1.1cm}Set Isnotgra:=Isnotgra+1 and go to Step 4; otherwise, compute the search direction ${d_{k + 1}}$ by (\ref{eq:3.6}) with (\ref{eq:3.8}) and \hspace*{1.1cm}(\ref{eq:3.9}). Set Isnotgra:=Isnotgra+1 and go to Step 4. \\
	 \hspace*{0.4cm} \textbf{3.3.} If the condition (\ref{eq:3.23}) holds, compute the search direction ${d_{k + 1}}$ by (\ref{eq:1.3}) where ${\beta _k} = \beta _k^{HS}.$ Set Isnotgra:=Isnotgra+1 \hspace*{1.1cm}and go to Step 4; otherwise, compute the search direction ${d_{k + 1}}$ by (\ref{eq:3.24}). Set Numgrad := Numgrad+1, Isnotgra \hspace*{1.1cm}:=0 and IterRestart :=0, and go to Step 4.\\              
	\textbf{Step 4.} Update ${Q_{k + 1}}$ and ${C_{k + 1}}$ using (\ref{eq:1.7}) and (\ref{eq:1.6}) with (\ref{eq:1.8}).\\
	\textbf{Step 5.} Set $k: = k + 1,$ and go to Step 1.\\
\end{algorithm}
\noindent\textbf{Remark 3} In Algorithm 1, Numgrad denotes the number of the successive use of the negative gradient direction; 
Isnotgra denotes the number of the successive use of the CG direction;
MaxRestart represents a quantification and when
 the Isnotgra reaches this value, we restart the method with $ - {g_k}$;
MinQuad also represents a quantification and when the IterQuad reaches this value, we restart the method with $ - {g_k}$. These parameters are related to the restart of the algorithm, which has an important impact on the numerical performance of the CG.

Secondly, we describe the subspace minimization conjugate gradient method in which the direction of the regularization model is generated by the problem (\ref{eq:3.14}).

If the condition
\begin{equation}\label{eq:4.2}
{(g_k^T{s_{k - 1}})^2} > (1 - {10^{ - 5}}){\left\| {{g_k}} \right\|^2}{\left\| {{s_{k - 1}}} \right\|^2}
\end{equation}
holds, the value of $\frac{{{{\left( {g_k^T{s_{k - 1}}} \right)}^2}}}{{{{\left\| {{g_k}} \right\|}^2}{{\left\| {{s_{k - 1}}} \right\|}^2}}}$ is close to 1, which means that vectors ${g_k}$ and ${s_{k - 1}}$ may be linearly correlated. So the positive definiteness of the matrix ${E_k}$ in (\ref{eq:3.14}) might not be guaranteed. Therefore, we choose the quadratic model to derive a search direction. 

We may consider to use ``\textbf{3.2.} If the condition (\ref{eq:3.2}) or (\ref{eq:3.4}) or (\ref{eq:4.1}) holds, compute the search direction ${d_{k + 1}}$ by (\ref{eq:3.6}) with (\ref{eq:3.21}) and (\ref{eq:3.22}). Set Isnotgra:=Isnotgra+1 and go to Step 4; otherwise, if the condition (\ref{eq:4.2}) holds, compute the search direction ${d_{k + 1}}$ by (\ref{eq:3.6}), (\ref{eq:3.15}) and (\ref{eq:3.16}) with $\lambda  = 0$, otherwise, compute the search direction ${d_{k + 1}}$ by (\ref{eq:3.6}) with (\ref{eq:3.15}) and (\ref{eq:3.16}). Set Isnotgra:=Isnotgra+1 and go to Step 4." to replace the Step 3.2 in Algorithm 1. The resulting method is called SMCG\_PR2.  We use SMCG\_PR to denote either SMCG\_PR1 or SMCG\_PR2.
  
The following two Lemmas show some properties of the direction ${d_k},$ which are essential to the convergence of SMCG\_PR.\\
\textbf{Lemma 4.1} Suppose the direction ${d_k}$ is calculated by SMCG\_PR. Then, there exists a constant ${c_1}$ such that 
\begin{equation}\label{eq:4.3}
g_k^T{d_k} \le  - {c_1}{\left\| {{g_k}} \right\|^2}.
\end{equation}

\emph{Proof.} We divide the proof into four cases.

\textbf{Case 1.} The direction ${d_k}$ is given by (\ref{eq:3.6}) with (\ref{eq:3.8}) and (\ref{eq:3.9}), as in SMCG\_PR1. Denote $T = \frac{1}{{1 + {\sigma _k}{{\left( {{z^*}} \right)}^{p - 2}}}}.$ Obviously, in this case, 
\begin{align*}
{\mu _k} = T{{\bar \mu }_k},\;\;{\nu _k} = T{{\bar \nu }_k}.
\end{align*}
If  ${\sigma _k}{\left( {{z^*}} \right)^{p - 2}} > 1$, we have $T = \frac{1}{2}$ from the first line after(\ref{eq:3.13}). Moreover, ${\sigma _k}\left( {{z^*}} \right) \ge 0.$ So we can establish that $\frac{1}{2} \le T \le 1$. From (3.31) and (3.32) of \cite{18.}, we can get that
\begin{equation}\label{eq:5.09}
g_k^T{d_k} = Tg_k^T({{\bar \mu }_k} {g_k} + {{\bar \nu }_k}{s_{k - 1}}) \le  - T\frac{{{{\left\| {{g_k}} \right\|}^4}}}{{{\rho _k}}} \le  - \frac{{{{\left\| {{g_k}} \right\|}^4}}}{{2{\rho _k}}}.
\end{equation}
Substituting ${\rho _k} = \frac{3}{2}\frac{{{{\left\| {{y_{k - 1}}} \right\|}^2}}}{{s_{k - 1}^T{y_{k - 1}}}}{\left\| {{g_k}} \right\|^2}$ into (\ref{eq:5.09}), we deduce that 
$g_k^T{d_k} \le  - \frac{{{{\left\| {{g_k}} \right\|}^4}}}{{2{\rho _k}}} =  - \frac{1}{3}\frac{{s_{k - 1}^T{y_{k - 1}}}}{{{{\left\| {{y_{k - 1}}} \right\|}^2}}}{\left\| {{g_k}} \right\|^2}.$
From (\ref{eq:3.18}), we konw $ - \frac{1}{{{\xi _1}}} \le  - \frac{{s_{k - 1}^T{y_{k - 1}}}}{{{{\left\| {{y_{k - 1}}} \right\|}^2}}} \le  - \frac{1}{{{\xi _2}}}.$ Therefore, we get
\begin{equation}\label{eq:4.5}
g_k^T{d_k} \le  - \frac{{{{\left\| {{g_k}} \right\|}^4}}}{{2{\rho _k}}} =  - \frac{1}{3}\frac{{s_{k - 1}^T{y_{k - 1}}}}{{{{\left\| {{y_{k - 1}}} \right\|}^2}}}{\left\| {{g_k}} \right\|^2} \le  - \frac{1}{{3{\xi _2}}}{\left\| {{g_k}} \right\|^2}.
	\end{equation}\\
On the other hand, if the direction ${d_k}$ is given by (\ref{eq:3.6}) with (\ref{eq:3.15}) and (\ref{eq:3.16}), which in SMCG\_PR2. We have that by direct calculation
\begin{align*}
\begin{array}{l}
g_k^T{d_k} = {{\hat \mu }_k}{\left\| {{g_k}} \right\|^2} + {{\hat \nu }_k}g_k^T{s_{k - 1}}\\
\;\;\;\;\;\;\;\;{\rm{        }} =  - \frac{{{{\left\| {{g_k}} \right\|}^4}}}{{\overline {{\Delta _k}} }}\left( {s_{k - 1}^T{y_{k - 1}} - 2g_k^T{y_{k - 1}}\frac{{g_k^T{s_{k - 1}}}}{{{{\left\| {{g_k}} \right\|}^2}}} + {\rho _k}{{\left( {\frac{{g_k^T{s_{k - 1}}}}{{{{\left\| {{g_k}} \right\|}^2}}}} \right)}^2} - \lambda g_k^T{s_{k - 1}}\frac{{g_k^T{s_{k - 1}}}}{{{{\left\| {{g_k}} \right\|}^2}}} + \lambda {{\left\| {{s_{k - 1}}} \right\|}^2}} \right)\\
\;\;\;\;\;\;\;\;{\rm{        }} =  - \frac{{{{\left\| {{g_k}} \right\|}^4}}}{{\overline {{\Delta _k}} }}\left( {\left( {{\rho _k} + \lambda {{\left\| {{g_k}} \right\|}^2}} \right){{\left( {\frac{{g_k^T{s_{k - 1}}}}{{{{\left\| {{g_k}} \right\|}^2}}}} \right)}^2} - \left( {2g_k^T{y_{k - 1}} + 2\lambda g_k^T{s_{k - 1}}} \right)\frac{{g_k^T{s_{k - 1}}}}{{{{\left\| {{g_k}} \right\|}^2}}} + s_{k - 1}^T{y_{k - 1}} + \lambda {{\left\| {{s_{k - 1}}} \right\|}^2}} \right)\\
\;\;\;\;\;\;\;\;{\rm{        }} \le  - \frac{{{{\left\| {{g_k}} \right\|}^4}}}{{\overline {{\Delta _k}} }}\frac{{{{\bar \Delta }_k} }}{{{\rho _k} + \lambda {{\left\| {{g_k}} \right\|}^2}}}\\
\;\;\;\;\;\;\;\;{\rm{        }}
 = \frac{{ - {{\left\| {{g_k}} \right\|}^2}}}{{\frac{3}{2}\frac{{{{\left\| {{y_{k - 1}}} \right\|}^2}}}{{s_{k - 1}^T{y_{k - 1}}}} + \lambda }} \\
 \;\;\;\;\;\;\;\;\;\le  - \frac{2}{{5{\xi _2}}}{\left\| {{g_k}} \right\|^2}.
\end{array}
\end{align*}
Due to $0 \le \lambda  \le \frac{{{{\left\| {{y_{k - 1}}} \right\|}^2}}}{{s_{k - 1}^T{y_{k - 1}}}},$ we have $\frac{3}{2}\frac{{{{\left\| {{y_{k - 1}}} \right\|}^2}}}{{s_{k - 1}^T{y_{k - 1}}}} \le \frac{3}{2}\frac{{{{\left\| {{y_{k - 1}}} \right\|}^2}}}{{s_{k - 1}^T{y_{k - 1}}}} + \lambda  \le \frac{5}{2}\frac{{{{\left\| {{y_{k - 1}}} \right\|}^2}}}{{s_{k - 1}^T{y_{k - 1}}}}.$ So, $ - \frac{2}{3}\frac{{s_{k - 1}^T{y_{k - 1}}}}{{{{\left\| {{y_{k - 1}}} \right\|}^2}}} \le \frac{{ - 1}}{{\frac{3}{2}\frac{{{{\left\| {{y_{k - 1}}} \right\|}^2}}}{{s_{k - 1}^T{y_{k - 1}}}} + \lambda }} \le  - \frac{2}{5}\frac{{s_{k - 1}^T{y_{k - 1}}}}{{{{\left\| {{y_{k - 1}}} \right\|}^2}}}.$ From (\ref{eq:3.18}), we konw $ - \frac{1}{{{\xi _1}}} \le  - \frac{{s_{k - 1}^T{y_{k - 1}}}}{{{{\left\| {{y_{k - 1}}} \right\|}^2}}} \le  - \frac{1}{{{\xi _2}}}.$ Therefore, the last inequality is established.

\textbf{Case 2.} ${d_k} = {{\bar \mu }_k} {g_k} + {{\bar \nu }_k} {s_{k - 1}},$ where ${{\bar \mu }_k} $ and ${{\bar \nu }_k} $ are calculated by (\ref{eq:3.21}) and (\ref{eq:3.22}), respectively. From (\ref{eq:5.09}) and (\ref{eq:4.5}), we can get that
\begin{equation}\label{eq:4.4}
g_k^T{d_k} = g_k^T({{\bar \mu }_k} {g_k} + {{\bar \nu }_k} {s_{k - 1}}) \le  - \frac{2}{{3{\xi _2}}}{\left\| {{g_k}} \right\|^2}.
\end{equation}

\textbf{Case 3.} If the direction ${d_k}$ is given by (\ref{eq:1.3}) where ${\beta _k} = \beta _k^{HS},$ (\ref{eq:4.3}) is satisfied by setting ${c_1} = 1 - {\xi _3}.$ The proof is similar to Lemma 3 in \cite{33.}.

\textbf{Case 4.} As ${d_k} =  - {g_k},$ we can easily derive $g_k^T{d_k} =  - {\left\| {{g_k}} \right\|^2}$ which satisfies (\ref{eq:4.3}) by setting ${c_1} = \frac{1}{2}.$

To sum up, the sufficient descent condition (\ref{eq:4.3}) holds by setting
\begin{align*}
{c_1} = \min \left\{ {\frac{1}{2},1 - {\xi _3},\frac{2}{{3{\xi _2}}},\frac{1}{{3{\xi _2}}},\frac{2}{{5{\xi _2}}}} \right\},
\end{align*}
which completes the proof.\\
\textbf{Lemma 4.2} Suppose the direction ${d_k}$ is calculated by SMCG\_PR. Then, there exists a constant ${c_2} > 0$ such that 
\begin{equation}\label{eq:4.6}
\left\| {{d_k}} \right\| \le {c_2}\left\| {{g_k}} \right\|.
\end{equation}
\emph{Proof.} The proof is also divided into four parts.

\textbf{Case 1.} The direction ${d_k}$ is given by (\ref{eq:3.6}) with (\ref{eq:3.8}) and (\ref{eq:3.9}), as in SMCG\_PR1.
 From (3.12) in \cite{33.} and $T \le 1$, we obtain
 \begin{align*}
 \left\| {{d_k}} \right\| = T\left\| {{{\bar \mu }_k} {g_k} + {{\bar \nu }_k} {s_{k - 1}}} \right\| \le \frac{{20}}{{{\xi _1}}}\left\| {{g_k}} \right\|.
 \end{align*}

On the other hand, if the direction ${d_k}$ is given by (\ref{eq:3.6}) with (\ref{eq:3.15}) and (\ref{eq:3.16}), as in SMCG\_PR2. At first, we give a lower bound of ${{\bar \Delta }_k}.$ From (\ref{eq:3.17}), we have
	\begin{align*}
	{\bar \Delta _k} = {\lambda ^2}\left( {{{\left\| {{g_k}} \right\|}^2}{{\left\| {{s_{k - 1}}} \right\|}^2} - {{\left( {g_k^T{s_{k - 1}}} \right)}^2}} \right) + \lambda \left( {{\rho _k}{{\left\| {{s_{k - 1}}} \right\|}^2} + s_{k - 1}^T{y_{k - 1}}{{\left\| {{g_k}} \right\|}^2} - 2g_k^T{y_{k - 1}}g_k^T{s_{k - 1}}} \right) +\\ {\rho _k}s_{k - 1}^T{y_{k - 1}} - {\left( {g_k^T{y_{k - 1}}} \right)^2}.
	\end{align*} 
	Moreover, using the Cauchy inequality and average inequality, we have
	\begin{align*}
	\begin{array}{l}
	\;\;\;\;{\rho _k}{\left\| {{s_{k - 1}}} \right\|^2} + s_{k - 1}^T{y_{k - 1}}{\left\| {{g_k}} \right\|^2} - 2g_k^T{y_{k - 1}}g_k^T{s_{k - 1}}\\
	\ge \frac{3}{2}\frac{{{{\left\| {{y_{k - 1}}} \right\|}^2}{{\left\| {{s_{k - 1}}} \right\|}^2}}}{{s_{k - 1}^T{y_{k - 1}}}}{\left\| {{g_k}} \right\|^2} + s_{k - 1}^T{y_{k - 1}}{\left\| {{g_k}} \right\|^2} - 2\left\| {{s_{k - 1}}} \right\|\left\| {{y_{k - 1}}} \right\|{\left\| {{g_k}} \right\|^2}\\
	= \left( {\frac{1}{2}\frac{{\left\| {{s_{k - 1}}} \right\|\left\| {{y_{k - 1}}} \right\|}}{{s_{k - 1}^T{y_{k - 1}}}} + \frac{{\left\| {{s_{k - 1}}} \right\|\left\| {{y_{k - 1}}} \right\|}}{{s_{k - 1}^T{y_{k - 1}}}} + \frac{{s_{k - 1}^T{y_{k - 1}}}}{{\left\| {{s_{k - 1}}} \right\|\left\| {{y_{k - 1}}} \right\|}} - 2} \right)\left\| {{s_{k - 1}}} \right\|\left\| {{y_{k - 1}}} \right\|{\left\| {{g_k}} \right\|^2}\\
	\ge \left( {\frac{1}{2}\frac{{\left\| {{s_{k - 1}}} \right\|\left\| {{y_{k - 1}}} \right\|}}{{s_{k - 1}^T{y_{k - 1}}}} + 2 - 2} \right)\left\| {{s_{k - 1}}} \right\|\left\| {{y_{k - 1}}} \right\|{\left\| {{g_k}} \right\|^2}\\
	\ge \frac{1}{2}\left\| {{s_{k - 1}}} \right\|\left\| {{y_{k - 1}}} \right\|{\left\| {{g_k}} \right\|^2} \ge 0.
	\end{array}
	\end{align*}
	It follows from (\ref{eq:3.18}) that $s_{k - 1}^T{y_{k - 1}} \ge {\xi _1}{\left\| {{s_{k - 1}}} \right\|^2}.$ By ${\rho _k} = \frac{3}{2}\frac{{{{\left\| {{y_{k - 1}}} \right\|}^2}}}{{s_{k - 1}^T{y_{k - 1}}}}{\left\| {{g_k}} \right\|^2},\;\;\lambda  \ge 0$ and the Cauchy inequality, we obtain a lower bound of ${{\bar \Delta }_k}$ that
	\begin{align*}
	\begin{array}{l}
	{{\bar \Delta }_k} \ge {\rho _k}s_{k - 1}^T{y_{k - 1}} - {\left( {g_k^T{y_{k - 1}}} \right)^2} = s_{k - 1}^T{y_{k - 1}}\left( {{\rho _k} - \frac{{{{\left( {g_k^T{y_{k - 1}}} \right)}^2}}}{{s_{k - 1}^T{y_{k - 1}}}}} \right)\\
	\;\;\;\;\;\ge {\xi _1}{\left\| {{s_{k - 1}}} \right\|^2}\left( {{\rho _k} - \frac{{{{\left( {g_k^T{y_{k - 1}}} \right)}^2}}}{{s_{k - 1}^T{y_{k - 1}}}}} \right)\\
	\;\;\;\;\;\ge \frac{1}{2}{\xi _1}{\left\| {{s_{k - 1}}} \right\|^2}\frac{{{{\left\| {{y_{k - 1}}} \right\|}^2}}}{{s_{k - 1}^T{y_{k - 1}}}}{\left\| {{g_k}} \right\|^2}.
	\end{array}
	\end{align*}  
	Using the triangle inequality, Cauchy inequality, ${\rho _k} = \frac{3}{2}\frac{{{{\left\| {{y_{k - 1}}} \right\|}^2}}}{{s_{k - 1}^T{y_{k - 1}}}}{\left\| {{g_k}} \right\|^2},$ $0 \le \lambda  \le \frac{{{{\left\| {{y_{k - 1}}} \right\|}^2}}}{{s_{k - 1}^T{y_{k - 1}}}}$ and the last relation, we have
	\begin{align*}
	\begin{array}{l}
	\left\| {{d_k}} \right\| = \left\| {{{\hat \mu }_k}{g_k} + {{\hat \nu }_k}{s_{k - 1}}} \right\|\\
	= \left\| {\frac{1}{{{{\bar \Delta }_k}}}\left( {\left( {g_k^T{y_{k - 1}}g_k^T{s_{k - 1}} - s_{k - 1}^T{y_{k - 1}}{{\left\| {{g_k}} \right\|}^2} + \lambda \left( {{{\left( {g_k^T{s_{k - 1}}} \right)}^2} - {{\left\| {{s_{k - 1}}} \right\|}^2}{{\left\| {{g_k}} \right\|}^2}} \right)} \right){g_k} + \left( {g_k^T{y_{k - 1}}{{\left\| {{g_k}} \right\|}^2} - {\rho _k}g_k^T{s_{k - 1}}} \right){s_{k - 1}}} \right)} \right\|\\
	\le \frac{1}{{{{\bar \Delta }_k}}}\left( {\left( {\left| {g_k^T{y_{k - 1}}g_k^T{s_{k - 1}}} \right| + \left| {s_{k - 1}^T{y_{k - 1}}} \right|{{\left\| {{g_k}} \right\|}^2} + \lambda {{\left| {g_k^T{s_{k - 1}}} \right|}^2} + \lambda {{\left\| {{s_{k - 1}}} \right\|}^2}{{\left\| {{g_k}} \right\|}^2}} \right)\left\| {{g_k}} \right\| + \left| {g_k^T{y_{k - 1}}{{\left\| {{g_k}} \right\|}^2} - {\rho _k}g_k^T{s_{k - 1}}} \right|\left\| {{s_{k - 1}}} \right\|} \right)\\
	\le \frac{1}{{{{\bar \Delta }_k}}}\left( {\left( {2\left\| {{s_{k - 1}}} \right\|\left\| {{y_{k - 1}}} \right\| + 2\frac{{{{\left\| {{y_{k - 1}}} \right\|}^2}{{\left\| {{s_{k - 1}}} \right\|}^2}}}{{s_{k - 1}^T{y_{k - 1}}}}} \right){{\left\| {{g_k}} \right\|}^3} + \left( {\left\| {{s_{k - 1}}} \right\|\left\| {{y_{k - 1}}} \right\| + \frac{{{\rho _k}}}{{{{\left\| {{g_k}} \right\|}^2}}}{{\left\| {{s_{k - 1}}} \right\|}^2}} \right){{\left\| {{g_k}} \right\|}^3}} \right)\\
	= \frac{1}{{{{\bar \Delta }_k}}}\left( {\left( {3\left\| {{s_{k - 1}}} \right\|\left\| {{y_{k - 1}}} \right\| + \frac{7}{2}\frac{{{{\left\| {{y_{k - 1}}} \right\|}^2}{{\left\| {{s_{k - 1}}} \right\|}^2}}}{{s_{k - 1}^T{y_{k - 1}}}}} \right){{\left\| {{g_k}} \right\|}^3}} \right)\\
	\le \frac{{13}}{{{\xi _1}}}\left\| {{g_k}} \right\|.
	\end{array}
	\end{align*}

\textbf{Case 2.} ${d_k} = {{\bar \mu }_k} {g_k} + {{\bar \nu }_k} {s_{k - 1}},$ where ${{\bar \mu }_k} $ and ${{\bar \nu }_k} $ are calculated by (\ref{eq:3.21}) and (\ref{eq:3.22}), respectively. From (3.12) in \cite{33.}, we can get (\ref{eq:4.6}) is satisfied by setting ${c_2} = \frac{{20}}{{{\xi _1}}}.$

\textbf{Case 3.} If the direction ${d_k}$ is given by (\ref{eq:1.3}) where ${\beta _k} = \beta _k^{HS},$ (\ref{eq:4.6}) is satisfied by setting ${c_2} = 1 + \frac{L}{{{\xi _1}}}.$ The proof is same as Lemma 4 in \cite{33.}.

\textbf{Case 4.} As ${d_k} =  - {g_k},$ we can easily establish that $\left\| {{d_k}} \right\| = \left\| {{g_k}} \right\|.$

In summary, we easily obtain the fact that (\ref{eq:4.6}) holds by 
\begin{align*}
{c_2} = \max \left\{ {1,1 + \frac{L}{{{\xi _1}}},\frac{{20}}{{{\xi _1}}}} \right\},
\end{align*}
which completes the proof.
\section{Convergence Analysis}
\label{sec:5}
In this section, we establish the global convergence and $R - $linear convergence of SMCG\_PR. We assume that $\left\| {{g_k}} \right\| \ne 0$ for each $k;$ otherwise, there is a stationary point for some $k.$

At first, we suppose that the objective function $f$ satisfies the following assumptions. Define $\Theta $ as an open neighborhood of the level set $L\left( {{x_0}} \right) = \left\{ {x \in {R^n}:f\left( x \right) \le f\left( {{x_0}} \right)} \right\},$ where ${{x_0}}$ is the initial point.\\
\textbf{Assumption 1} $f$ is continuously differentiable and bounded from below in $\Theta. $\\
\textbf{Assumption 2} The gradient $g$ is Lipchitz continuous in $\Theta ,$ namely, there exists a constant $L > 0$ such that $\left\| {g(x) - g(y)} \right\| \le L\left\| {x - y} \right\|,\forall x,y \in \Theta .$\\
\textbf{Lemma 5.1} Suppose the Assumption 1 holds and  the iterative sequence $\{ {x_k}\} $ is generated by the SMCG\_PR. Then, we have ${f_k} \le {C_k}$ for each $k.$\\
\emph{Proof.} Due to (\ref{eq:1.4}) and descent direction ${d_{k + 1}},$ ${f_{k + 1}} < {C_k}$ always holds. Through (\ref{eq:1.7}), we can get ${C_1} = {C_0}$ or ${C_1} = {f_1} + 1.0.$ If ${C_1} = {C_0},$ because of the relations ${f_{k + 1}} < {C_k}$  and ${C_0} = {f_0},$ we know ${f_1} \le {C_1}.$ If ${C_1} = {f_1} + 1.0,$ we can easily get ${f_1} \le {C_1}.$ When $k \ge 1,$ the updated form of ${C_{k + 1}}$ is (\ref{eq:1.6}), similar to Lemma 1.1 in \cite{52.}, we have ${f_{k + 1}} \le {C_{k + 1}}.$ Therefore, ${f_k} \le {C_k}$ holds for each $k.$ \\
\textbf{Lemma 5.2} Suppose the Assumption 2 holds and the iterative sequence $\{ {x_k}\} $ is generated by the SMCG\_PR. Then, 
\begin{equation}\label{eq:5.1}
{\alpha _k} \ge \left( {\frac{{1 - \sigma }}{L}} \right)\frac{{\left| {g_k^T{d_k}} \right|}}{{{{\left\| {{d_k}} \right\|}^2}}}.
\end{equation}
\emph{Proof.} By (\ref{eq:1.5}) and Assumption 2, we have that
	\begin{align*}
	\left( {\sigma  - 1} \right)g_k^T{d_k} \le {\left( {{g_{k + 1}} - {g_k}} \right)^T}{d_k} \le {\alpha _k}L{\left\| {{d_k}} \right\|^2}.
	\end{align*}
	Since ${d_k}$ is a descent direction and $\sigma  < 1,$ (\ref{eq:5.1}) follows immediately.\\
\textbf{Theorem 5.3} Suppose Assumption 1 and 2 hold. If the iterative sequence $\{ {x_k}\} $ is generated by the SMCG\_PR, it follows
\begin{equation}\label{eq:5.2}
\mathop {\lim }\limits_{k \to \infty } \left\| {g({x_k})} \right\| = 0.
\end{equation}
\emph{Proof.} By (\ref{eq:1.4}), Lemma 5.2, Lemma 4.1, and Lemma 4.2, we get that
\begin{align*}
{f_{k + 1}} \le {C_k} - \frac{{\delta \left( {1 - \sigma } \right)}}{L}\frac{{{{\left( {g_k^T{d_k}} \right)}^2}}}{{{{\left\| {{d_k}} \right\|}^2}}} \le {C_k} - \frac{{\delta \left( {1 - \sigma } \right)c_1^2}}{{Lc_2^2}}{\left\| {{g_k}} \right\|^2}.
\end{align*}
In short, set $\beta  = \frac{{\delta \left( {1 - \sigma } \right)c_1^2}}{{Lc_2^2}},$ we give the fact that
\begin{equation}\label{eq:5.3}
{f_{k + 1}} \le {C_k} - \beta {\left\| {{g_k}} \right\|^2}.
\end{equation}

Now, we find a upper bound of ${Q_{k + 1}}$ in (\ref{eq:1.6}) with (\ref{eq:1.8}). As for $k \ge 1,$ ${Q_{k + 1}}$ can be expressed as \textcolor{red}{\cite{35.}}
\begin{align*}
{Q_{k + 1}} = \left\{ \begin{array}{l}
1 + (l + 1)\sum\limits_{i = 1}^{k/l} {{\eta ^i}} ,{\rm{\;\;\;\;\;\;\;\;\;\;\;\;\;\;\;\;\;\;\;\;\;}}\bmod (k,l) = 0,\\
1 + \bmod (k,l) + (l + 1)\sum\limits_{i = 1}^{\left\lfloor {k/l} \right\rfloor } {{\eta ^i}} ,{\rm{ }}\bmod (k,l) \ne 0,
\end{array} \right.
\end{align*}
where $\left\lfloor . \right\rfloor $ is the floor function. Then, we obtain
\begin{align}\label{eq:5.4}
\begin{array}{l}
{Q_{k + 1}} \le 1 + \bmod (k,l) + (l + 1)\sum\limits_{i = 1}^{\left\lfloor {k/l} \right\rfloor  + 1} {{\eta ^i}} \\
\;\;\;\;\;\;\;\;\;\;\le 1 + (l + 1) + (l + 1)\sum\limits_{i = 1}^{\left\lfloor {k/l} \right\rfloor  + 1} {{\eta ^i}} \\
\;\;\;\;\;\;\;\;\;\;\le 1 + \left( {l + 1} \right) + \left( {l + 1} \right)\sum\limits_{i = 1}^{k + 1} {{\eta ^i}}  \\
\;\;\;\;\;\;\;\;\;\;= 1 + (l + 1)\sum\limits_{i = 0}^{k + 1} {{\eta ^i}} \\
\;\;\;\;\;\;\;\;\;\;= 1 + \frac{{(l + 1)(1 - {\eta ^{k + 2}})}}{{1 - \eta }}\\
\;\;\;\;\;\;\;\;\;\;\le 1 + \frac{{l + 1}}{{1 - \eta }}.
\end{array}
\end{align}
Denote $M = 1 + \frac{{l + 1}}{{1 - \eta }},$ which gives the fact ${Q_{k + 1}} \le M.$\\
With the updated form of ${C_{k + 1}}$ in (\ref{eq:1.6}), (\ref{eq:5.3}) and (\ref{eq:5.4}), we obtain
\begin{equation}\label{eq:80}
{C_{k + 1}} = {C_k} + \frac{{{f_{k + 1}} - {C_k}}}{{{Q_{k + 1}}}} \le {C_k} - \frac{\beta }{{{Q_{k + 1}}}}{\left\| {{g_k}} \right\|^2} \le {C_k} - \frac{\beta }{M}{\left\| {{g_k}} \right\|^2}.
\end{equation}
According to (\ref{eq:1.7}), we know ${C_1} \le {C_0}$ which implies that ${C_k}$ is monotonically decreasing. Due to Assumption 1 and Lemma 5.1, we can get ${C_k}$ is bounded from below. Then
\begin{align*}
\sum\limits_{k = 0}^\infty  {\frac{\beta }{M}{{\left\| {{g_k}} \right\|}^2} < \infty } ,
\end{align*} 
therefore,
\begin{align*}
\mathop {\lim }\limits_{k \to \infty } \left\| {g({x_k})} \right\| = 0,
\end{align*}
which completes the proof.

Moreover, $R - $linear convergence of SMCG\_PR will be established as followed. In order to establish $R - $linear convergence of SMCG\_PR, we introduce Definition 1 and assume that the optimal set ${\chi ^*}$ is nonempty.\\
\textbf{Definition 1} The continuously differentiable function $f$ has a global error bound on ${R^n}$, if there exists a constant ${\kappa _f} > 0$ such that for any $x \in {R^n}$ and $\overline x  = {[x]_{{\chi ^*}}}$, we have 
\begin{equation}\label{eq:5.5}
\left\| {g\left( x \right)} \right\| \ge {\kappa _f}\left\| {x - \overline x } \right\|\;\;\forall x \in {R^n},
\end{equation}
where $\overline x  = {[x]_{{\chi ^*}}}$ is the projection of $x$ onto the nonempty solution set ${\chi ^*}.$ We further denote by ${\chi ^*} = \arg {\min _{x \in {R^n}}}f\left( x \right)$ the set of optimal solutions of problem (\ref{eq:1.1}).\\
	\textbf{Remark 4} By Assumption 2 it is $\left\| {g\left( x \right) - g\left( {{x^*}} \right)} \right\| \le L\left\| {x - {x^*}} \right\|,$ so that it is also $\left\| {g\left( x \right)} \right\| \le L\left\| {x - {x^*}} \right\|,$ which implies ${k_f} \le L.$\\
\textbf{Remark 5} \textcolor{red}{\cite{555.}} If $f$ is strongly convex, it must satisfy Definition 1.\\
\textbf{Remark 6} If $f$ is a convex function and the optimal solution set is nonempty, the function value at the optimal solution is equal.\\ 
\textbf{Theorem 5.4} Suppose that Assumption 2 holds, $f$ is convex with a minimizer ${x^*}$ and the solution set ${\chi ^*}$ is nonempty, and there exists $\overline \alpha   > 0$ such that ${\alpha _k} \le \overline \alpha  $ for all $k.$ Let $f$ satisfy Definition 1 with constant ${\kappa _f} > 0.$ In what follows, we only consider the case of $\left\| {{g_k}} \right\| \ne 0,$ $\forall k \ge 0.$ Then there exists $\theta  \in (0,1)$ such that 
\begin{align*}
{f_k} - f({x^*}) \le {\theta ^k}({f_0} - f({x^*})).
\end{align*}
\emph{Proof.} From Lemma 5.1, we can get ${f_{k + 1}} \le {C_{k + 1}}.$ Due to Remark 6 and $\left\| {{g_k}} \right\| \ne 0,$ $\forall k \ge 0$, we know ${x_{k + 1}}$ is not the optimal solution. So, we have $f\left( {{x^*}} \right) < {f_{k + 1}}.$ From (\ref{eq:80}) and $\left\| {{g_k}} \right\| \ne 0,$ $\forall k \ge 0,$ we have that ${C_{k + 1}} < {C_k}.$ Therefore, we get $f\left( {{x^*}} \right) < {f_{k + 1}} \le {C_{k + 1}} < {C_k},$ which means $f\left( {{x^*}} \right) < {C_{k + 1}} < {C_k}.$ It follows
\begin{equation}\label{eq:5.6}
0 < \frac{{{C_{k + 1}} - f({x^*})}}{{{C_k} - f({x^*})}} < 1,\;\;\forall k \ge 0.
\end{equation}
Set
\begin{equation}\label{eq:5.7}
r = \mathop {\lim }\limits_{k \to \infty } {\rm{ sup}}\frac{{{C_{k + 1}} - f({x^*})}}{{{C_k} - f({x^*})}},
\end{equation}
then, $0 \le r \le 1.$

First of all, we consider the case of $r = 1.$ According to (\ref{eq:5.7}), there exists a subsequence $\{ {x_{{k_j}}}\} $ such that
\begin{equation}\label{eq:5.8}
\mathop {\lim }\limits_{j \to \infty } \frac{{{C_{{k_j} + 1}} - f({x^*})}}{{{C_{{k_j}}} - f({x^*})}} = 1.
\end{equation} 

Because of (\ref{eq:5.4}), there exists $q > 0,$ $0 < q \le \frac{1}{{{Q_{{k_j} + 1}}}} \le 1$ holds. Hence, there exists a subsequence of $\{ {x_{{k_j}}}\} $ such that the corresponding subsequence of $\left\{ {\frac{1}{{{Q_{{k_j} + 1}}}}} \right\}$ is convergent. Without loss of generality, we assume that 
\begin{equation}\label{eq:5.9}
\mathop {\lim }\limits_{j \to \infty } \frac{1}{{{Q_{{k_j} + 1}}}} = {r_1}.
\end{equation}
Clearly, $0 < {r_1} \le 1.$

By the updating formula of ${C_{k + 1}}$ in (\ref{eq:1.6}), we obtain
\begin{align*}
\frac{{{C_{{k_j} + 1}} - f({x^*})}}{{{C_{{k_j}}} - f({x^*})}} = \left( {1 - \frac{1}{{{Q_{{k_j} + 1}}}}} \right) + \frac{1}{{{Q_{{k_j} + 1}}}}\frac{{{f_{{k_j} + 1}} - f({x^*})}}{{{C_{{k_j}}} - f({x^*})}}.
\end{align*}
It follows from (\ref{eq:5.8}), (\ref{eq:5.9}) and finding the limit of upper formula that
\begin{equation}\label{eq:5.10}
\mathop {\lim }\limits_{j \to \infty } \frac{{{f_{{k_j} + 1}} - f({x^*})}}{{{C_{{k_j}}} - f({x^*})}} = 1.
\end{equation}
Using convexity of $f,$ the solution set ${\chi ^*}$ is nonempty and Remark 6, we know $f\left( {{x^*}} \right) = f\left( {\bar x} \right),$ where ${\bar x}$ is introduced in Definition 1. So, we have that ${f_{{k_j} + 1}} - f\left( {{x^*}} \right) = {f_{{k_j} + 1}} - f\left( {\bar x} \right).$ Through convexity of $f,$ we have ${f_{{k_j} + 1}} - f\left( {\bar x} \right) \le \left( {\nabla {f_{{k_j} + 1}},{x_{{k_j} + 1}} - \bar x} \right).$ According to Definition 1 and Cauchy-Schwarz inequality, then $\left( {\nabla {f_{{k_j} + 1}},{x_{{k_j} + 1}} - \bar x} \right) \le \frac{1}{{{k_f}}}{\left\| {{g_{{k_j} + 1}}} \right\|^2}.$ Therefore, we get
\begin{equation}\label{eq:5.11}
{f_{{k_j} + 1}} - f({x^*}) = {f_{{k_j} + 1}} - f(\overline x ) \le (\nabla {f_{{k_j} + 1}},{x_{{k_j} + 1}} - \overline x ) \le \frac{1}{{{\kappa _f}}}{\left\| {{g_{{k_j} + 1}}} \right\|^2}.
\end{equation}
According to the Lipschitz continuity of $g,$ ${\alpha _k} \le \overline \alpha  $ and (\ref{eq:4.6}), we have
\begin{align*}
\left\| {{g_{{k_j} + 1}}} \right\| \le \left\| {{g_{{k_j} + 1}} - {g_{{k_j}}}} \right\| + \left\| {{g_{{k_j}}}} \right\| \le L\left\| {{x_{{k_j} + 1}} - {x_{{k_j}}}} \right\| + \left\| {{g_{{k_j}}}} \right\| \le (1 + L\overline \alpha  {c_2})\left\| {{g_{{k_j}}}} \right\|,
\end{align*}
together with (\ref{eq:5.11}), it implies that
\begin{align*}
{f_{{k_j} + 1}} - f({x^*}) \le \frac{1}{{{\kappa _f}}}{(1 + L\overline \alpha  {c_2})^2}{\left\| {{g_{{k_j}}}} \right\|^2}.
\end{align*}
Dividing the above inequality by ${C_{{k_j}}} - f({x^*}),$ we have 
\begin{align}\label{eq:5.12}
0 < \frac{{{f_{{k_j} + 1}} - f({x^*})}}{{{C_{{k_j}}} - f({x^*})}} \le \frac{{{{(1 + L\overline \alpha  {c_2})}^2}{{\left\| {{g_{{k_j}}}} \right\|}^2}}}{{{\kappa _f}({C_{{k_j}}} - f({x^*}))}}.
\end{align}
Based on (\ref{eq:5.3})
\begin{align*}
{f_{{k_j} + 1}} - f({x^*}) \le {C_{{k_j}}} - f({x^*}) - \beta {\left\| {{g_{{k_j}}}} \right\|^2}.
\end{align*}
Dividing both sides of above inequality by ${C_{{k_j}}} - f({x^*}),$ we get
\begin{align*}
\frac{{{f_{{k_j} + 1}} - f({x^*})}}{{{C_{{k_j}}} - f({x^*})}} \le 1 - \frac{{\beta {{\left\| {{g_{{k_j}}}} \right\|}^2}}}{{{C_{{k_j}}} - f({x^*})}}.
\end{align*}
Combining with (\ref{eq:5.10}), then
\begin{align*}
\mathop {\lim }\limits_{j \to \infty } \frac{{{{\left\| {{g_{{k_j}}}} \right\|}^2}}}{{{C_{{k_j}}} - f({x^*})}} = 0,
\end{align*}
due to (\ref{eq:5.12}), it follows
\begin{align*}
\mathop {\lim }\limits_{j \to \infty } \frac{{{f_{{k_j} + 1}} - f({x^*})}}{{{C_{{k_j}}} - f({x^*})}} = 0,
\end{align*}
which contradicts with (\ref{eq:5.10}). Therefore, the case of $r = 1$ does not occur, that is,
\begin{align*}
\mathop {\lim }\limits_{k \to \infty } {\rm{ sup}}\frac{{{C_{k + 1}} - f({x^*})}}{{{C_k} - f({x^*})}} = r < 1.
\end{align*}
Then, there exists an integer ${k_0} > 0$ such that
\begin{equation}\label{eq:5.13}
\frac{{{C_{k + 1}} - f({x^*})}}{{{C_k} - f({x^*})}} \le r + \frac{{1 - r}}{2} = \frac{{1 + r}}{2} < 1,\;\;\forall k > {k_0}.
\end{equation}
From (\ref{eq:5.6}), we know that $0 < \mathop {\max }\limits_{0 \le k \le {k_0}} \left\{ {\frac{{{C_{k + 1}} - f({x^*})}}{{{C_k} - f({x^*})}}} \right\} = \overline r  < 1.$ Let $\theta  = max\left\{ {\frac{{1 + r}}{2},\overline r } \right\}.$\\
Clearly, $0 < \theta  < 1.$ It follows from (\ref{eq:5.13}) that
{\setlength\abovedisplayskip{1pt}
	\setlength\belowdisplayskip{1pt}
\begin{align*}
{C_{k + 1}} - f({x^*}) \le \theta \left( {{C_k} - f({x^*})} \right),
\end{align*}}
which indicates that
{\setlength\abovedisplayskip{1pt}
	\setlength\belowdisplayskip{1pt}
\begin{align*}
{C_{k + 1}} - f({x^*}) \le \theta \left( {{C_k} - f({x^*})} \right) \le {\theta ^{k + 1}}\left( {{C_0} - f({x^*})} \right).
\end{align*}}
In addition, due to ${f_{k + 1}} \le {C_{k + 1}}$ in Lemma 5.1 and ${C_0} = {f_0},$ we can deduce that
{\setlength\abovedisplayskip{1pt}
	\setlength\belowdisplayskip{1pt}
\begin{align*}
\left( {{f_k} - f({x^*})} \right) \le {\theta ^{k}}\left( {{f_0} - f({x^*})} \right),
\end{align*}}
which completes the proof.
\section{Numerical Results}
\label{sec:6}
In this section, numerical experiments are conducted to show the efficiency of the SMCG\_PR with $p = 3$ and $p = 4.$ We compare the performance of SMCG\_PR to that of CG\_DESCENT (5.3) \cite{28.}, CGOPT \cite{55.}, SMCG\_BB \cite{34.} and SMCG\_Conic \cite{54.} for the 145 test problems in the CUTEr library \cite{23.}. The names and dimensions for the 145 test problems are the same as that of the numerical results in \cite{30.}. The codes of CG\_DESCENT (5.3), CGOPT and SMCG\_BB  can be downloaded from \textcolor{blue}{http://users.clas.ufl.edu/hager/papers/Software}, \textcolor{blue}{http://coa.amss.ac.cn/wordpress/?page\_id=21} and \textcolor{blue}{http://web.xidian.edu.cn/xdliuhongwei/paper.html}, respectively.

The following parameters are used in SMCG\_PR:
\begin{align*}
\varepsilon  = {10^{ - 6}},\delta  = 0.0005,\sigma  = 0.9999,{\lambda _{\min }} = {10^{ - 30}},{\lambda _{\max }} = {10^{30}},\gamma  = {10^{ - 5}},
\end{align*}
\begin{align*}
{\xi _1} = {10^{ - 7}},{\xi _2} = 1.25 \times {10^4},{\xi _3} = {10^{ - 5}},{\xi _4} = {10^{ - 9}},{\xi _5} = {10^{ - 11}},{c_1} = {10^{ - 4}},{c_2} = 0.080.
\end{align*}
CG\_DESCENT (5.3), CGOPT, SMCG\_BB and SMCG\_Conic use the default parameters in their codes. All test methods are terminated if ${\left\| {{g_k}} \right\|_\infty } \le {10^{ - 6}}$ is satisfied or the number of iterations exceeds 200,000.
 
\setlength{\parskip}{0.0em} 
The performance profiles introduced by Dolan and Mor$\acute{\text{e}}$ \cite{19.} are used to display the performances of the test methods. We present three groups of the numerical experiments. They all run  in Ubuntu 10.04 LTS which is fixed in a VMware Workstation 10.0 installed in Windows 7. In the following Figs. 1-12 and Table 2, ``${N_{iter}}$",``${N_f}$",``${N_g}$" and ``${T_{cpu}}$" represent the number of iterations, the number of function evaluations, the number of gradient evaluations and CPU time(s), respectively.

In the first group of numerical experiments, we compare SMCG\_PR1 and SMCG\_PR2 with $p = 3$ and $p = 4$. All these test methods can successfully solve 139 problems. It is observed from Fig.1-Fig.4 that the SMCG\_PR1 with $p = 3$ is better than others.
\vspace{-0.2cm}
\begin{figure}[H]
	\centering
	\begin{minipage}[t]{0.49 \linewidth}
		\includegraphics[scale=0.50]{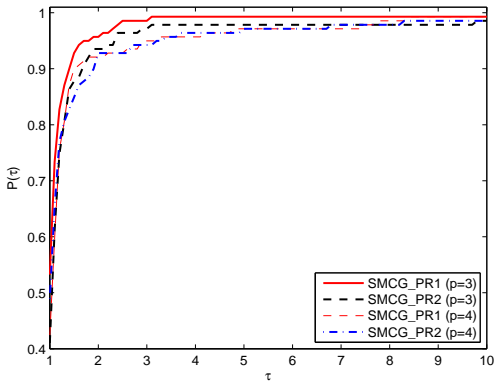}
		\setlength{\abovecaptionskip}{-0.2cm}
		\caption{Performance profile based on ${N_{iter}}$(CUTEr).}\label{fig.1}
	\end{minipage}	
	\begin{minipage}[t]{0.49\linewidth}
		\centering
		\includegraphics[scale=0.50]{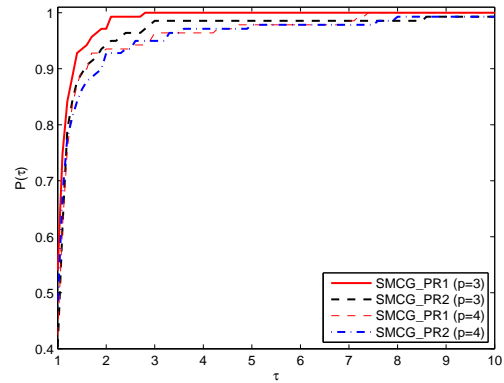}
		\setlength{\abovecaptionskip}{-0.2cm}
		\caption{Performance profile based on ${N_{f}}$(CUTEr).}\label{fig.2}
	\end{minipage}	
\end{figure}
\begin{figure}[H]
	\centering
	\begin{minipage}[t]{0.49 \linewidth}
		\includegraphics[scale=0.50]{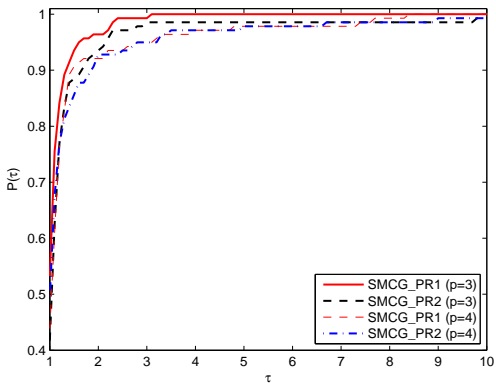}
		\setlength{\abovecaptionskip}{-0.2cm}
		\caption{Performance profile based on ${N_{g}}$(CUTEr).}\label{fig.3}
	\end{minipage}	
	\begin{minipage}[t]{0.49\linewidth}
		\centering
		\includegraphics[scale=0.50]{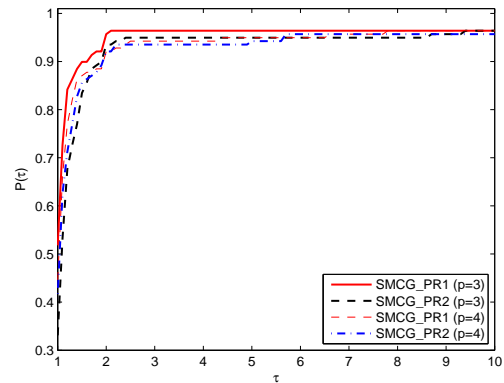}
		\setlength{\abovecaptionskip}{-0.2cm}
		\caption{Performance profile based on ${T_{cpu}}$(CUTEr).}\label{fig.4}
	\end{minipage}	
\end{figure}
In the second group of numerical experiments, we compare SMCG\_PR1 $\left( {p = 3} \right)$ with CG\_DECENT (5.3) and CGOPT. SMCG\_PR1 successfully solves 139 problems, while CG\_DECENT (5.3) and CGOPT successfully solve 144 and 134 problems, respectively.
\begin{figure}[H]
	\centering
	\begin{minipage}[t]{0.49 \linewidth}
		\includegraphics[scale=0.50]{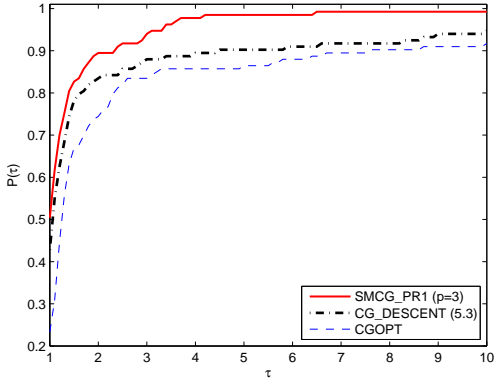}
		\caption{Performance profile based on ${N_{iter}}$(CUTEr).}\label{fig.5}
	\end{minipage}	
	\begin{minipage}[t]{0.49\linewidth}
		\centering
		\includegraphics[scale=0.50]{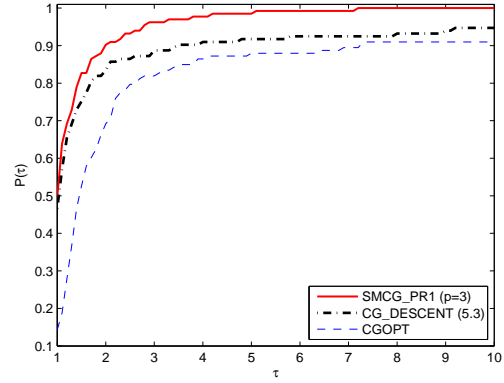}
		\caption{Performance profile based on ${N_{f}}$(CUTEr).}\label{fig.6}
	\end{minipage}	
\end{figure}

Regarding the number of iterations in Fig.5, we observe that SMCG\_PR1 is more efficient than CG\_DESCENT (5.3) and CGOPT, and it successfully solves about $50.4\% $ of the test problems with the least number of iterations, while the percentages of solved problems of CG\_DESCENT (5.3) and CGOPT are 42.8\% and 23.3\%, respectively. As shown in Fig.6, we see that SMCG\_PR1 outperforms CG\_DESCENT (5.3) and CGOPT for the number of function evaluations.
\begin{figure}[H]
	\centering
	\begin{minipage}[t]{0.49 \linewidth}
		\includegraphics[scale=0.50]{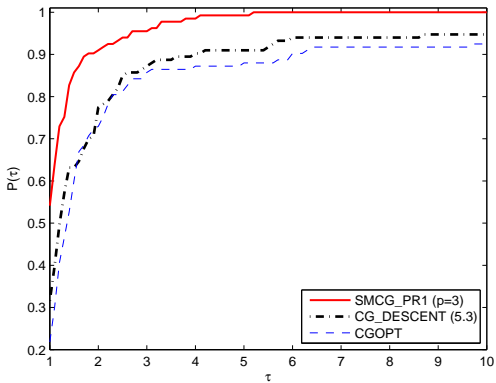}
		\caption{Performance profile based on ${N_{g}}$(CUTEr).}\label{fig.7}
	\end{minipage}	
	\begin{minipage}[t]{0.49\linewidth}
		\centering
		\includegraphics[scale=0.50]{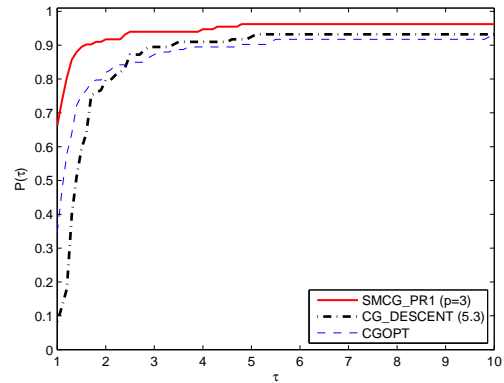}
		\caption{Performance profile based on ${T_{cpu}}$(CUTEr).}\label{fig.8}
	\end{minipage}	
\end{figure}

Fig.7 presents the performance profile relative to the number of gradient evaluations. We can observe that the SMCG\_PR1 is the top performance and solves about 54.2\% of test problems with the least number of gradient evaluations, and CG\_DESCENT (5.3) solves about 31.6\% and CGOPT solves about 21.8\%. From Fig.8, we can see that SMCG\_PR1 is fastest for about  66.2\%   of test problems, while CG\_DESCENT (5.3) and CGOPT are fastest for about 8.3\%    and  34.6\%, respectively. From Figs. 5, 6, 7 and 8, it indicates that SMCG\_PR1 outperforms CG\_DESCENT (5.3) and CGOPT for the 145 test problems in the CUTEr library. 

In the third group of the numerical experiments, we compare SMCG\_PR1 $\left( {p = 3} \right)$ with SMCG\_BB and SMCG\_Conic \cite{54.}. SMCG\_PR1 successfully solves 139 problems, which are 1 problem more than SMCG\_Conic, while SMCG\_BB successfully solves 140 problems. As shown in Figs. 9, 10, 11 and 12, we can easily observe that SMCG\_PR1 is superior to SMCG\_BB and SMCG\_Conic for the 145 test problems in the CUTEr library.

\begin{figure}[H]
	\centering
	\begin{minipage}[t]{0.49 \linewidth}
		\includegraphics[scale=0.50]{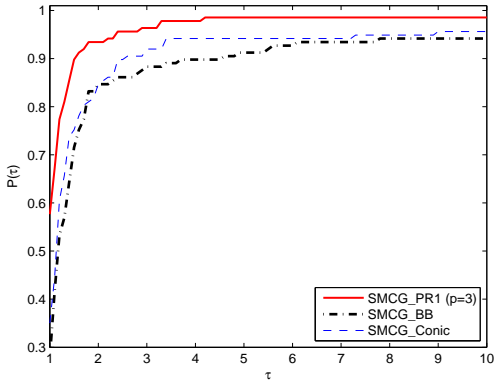}
		\caption{Performance profile based on ${N_{iter}}$(CUTEr).}\label{fig.9}
	\end{minipage}	
	\begin{minipage}[t]{0.49\linewidth}
		\centering
		\includegraphics[scale=0.50]{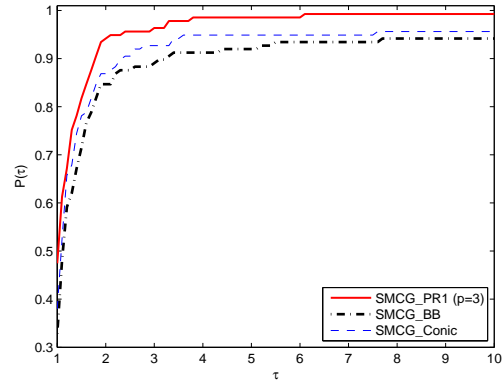}
		\caption{Performance profile based on ${N_{f}}$(CUTEr).}\label{fig.10}
	\end{minipage}	
\end{figure}
\begin{figure}[H]
	\centering
	\begin{minipage}[t]{0.49 \linewidth}
		\includegraphics[scale=0.50]{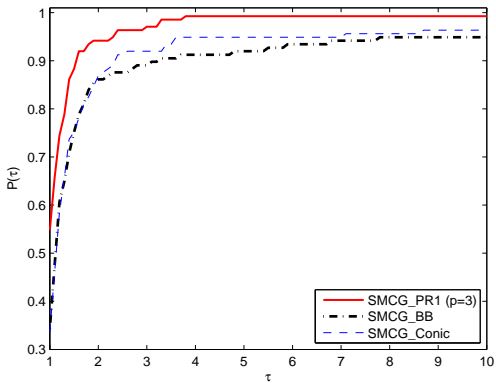}
		\caption{Performance profile based on ${N_{g}}$(CUTEr).}\label{fig.11}
	\end{minipage}	
	\begin{minipage}[t]{0.49\linewidth}
		\centering
		\includegraphics[scale=0.50]{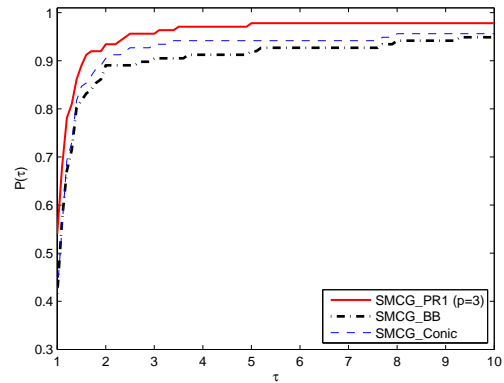}
		\caption{Performance profile based on ${T_{cpu}}$(CUTEr).}\label{fig.12}
	\end{minipage}	
\end{figure}

Due to limited space, we do not list all detailed numerical results. Instead, we present some numerical results about SMCG\_PR1 $\left( {p = 3} \right)$, CG\_DESCENT (5.3), CGOPT, SMCG\_BB and SMCG\_Conic for some ill-conditioned problems. Table 1 illustrates the notations, names and dimensions about the ill-conditioned problems. Table 2 presents some numerical results about  SMCG\_PR1 $\left( {p = 3} \right)$, CG\_DESCENT (5.3), CGOPT, SMCG\_BB and SMCG\_Conic for the problems in Table 1. As shown in Table 2, the most famous CG software packages CGOPT and CG\_DESCENT (5.3) both require many iterations, function evaluations and gradient evaluations when solving these ill-conditioned problems, though the dimensions of some of these  ill-conditioned problems are small. From Table 2, we observe that SMCG\_PR1 $\left( {p = 3} \right)$ has significant improvements over the other test methods, especially for CGOPT and CG\_DESCENT (5.3).  It indicates that SMCG\_PR1 $\left( {p = 3} \right)$ is  relatively competitive for ill-conditioned problems compared to other test methods.

\begin{table}[H]
	\caption{Some ill-conditioned problems in CUTEr}\centering
	\label{tab:1}       
	\begin{tabular}{cccccc}
		\hline\noalign{\smallskip}
		notation&name&dimension&notation&name&dimension  \\
		\noalign{\smallskip}\hline\noalign{\smallskip}
		P1&EIGENBLS&2550&P7&PALMER1D&7\\
		P2&EXTROSNB&1000&P8&PALMER2C&8\\
		P3&GROWTHLS&3&P9&PALMER4C&8\\
		P4&MARATOSB&2&P10&PALMER6C&8\\
		P5&NONCVXU2&5000&P11&PALMER7C&8\\
		P6&PALMER1C&8& & & \\	
		\noalign{\smallskip}\hline		 
	\end{tabular}
\end{table}
\begin{table}[H]
	\caption{Numerical results for some ill-conditioned problems in CUTEr}\centering
	\label{tab:2}       
	\begin{tabular}{cccccc}
		\hline\noalign{\smallskip}
		\multirow{2}{*}{problem}&SMCG\_PR1&CG\_DESCENT (5.3)&CGOPT&SMCG\_BB & SMCG\_Conic  \\
		\cline{2-6}
		& ${N_{iter}}/{N_f}/{N_g}$ & ${N_{iter}}/{N_f}/{N_g}$ & ${N_{iter}}/{N_f}/{N_g}$ & ${N_{iter}}/{N_f}/{N_g}$ &${N_{iter}}/{N_f}/{N_g}$\\
		\noalign{\smallskip}\hline\noalign{\smallskip}
		P1&9190/18382/9192&16092/32185/16093&19683/39369/19686&16040/32066/16041&12330/24654/12332\\
		P2&3568/6956/3574&6879/13839/6975&9127/18465/9305&8416/16195/8426&3733/7466/3735\\
		P3&1/2/2&441/997/596&480/1241/644&689/1512/711&1/2/2\\
		P4&212/614/389&946/2911/2191&1411/4185/2213&1159/9592/2634&3640/13621/5632\\
		P5&6096/12174/6098&7160/13436/8046&6195/12402/6207&6722/12800/6723&6459/12816/6460\\
		P6&1453/2093/1546&126827/224532/378489&Failed&88047/135548/89509&13007/23796/13352\\
		P7&445/682/470&3971/5428/10036&16490/36567/19846&2701/3703/2727&584/943/635     \\
		P8&307/440/318&21362/21455/42837&25716/61275/30492&4894/7169/5002&695/1386/697   \\
		P9&54/107/59&44211/49913/96429&88681/197232/105736& 1064/1622/1074&1055/2025/1071   \\
		P10&202/323/213&14174/14228/28411&29118/63118/31844&35704/58676/36281&1458/2429/1505   \\
		P11&6288/8757/6576&65294/78428/149585&98699/220388/119626&46397/65692/46929/&502/575/514  \\
		\noalign{\smallskip}\hline 
	\end{tabular}
\end{table}

The numerical results indicate that the SMCG\_PR method outperforms CG\_DESCENT (5.3), CGOPT, SMCG\_BB and SMCG\_Conic.
\section{Conclusions}
\label{sec:7}
In this paper, we present two new subspace minimization conjugate gradient methods based on the special $p - $regularization model for $p > 2.$ In the proposed methods, the search directions satisfy the sufficient descent condition. Under mild conditions, the global convergences of SMCG\_PR are established. We also prove that SMCG\_PR is $R - $linearly convergent. The numerical experiments show that SMCG\_PR is very promising.
\begin{acknowledgements}
We would like to thank  the anonymous referees for their useful  comments. We also would like to thank Professors Hager,W.W. and Zhang, H.C. for their CG\_DESCENT (5.3) code, and thank Professor Dai, Y.H and Dr. Kou, C.X. for their CGOPT code. This research is supported by National Science Foundation of China (No.11901561), Guangxi Natural Science Foundation (No.2018GXNSFBA281180) and China Postdoctoral Science Foundation (2019M660833).
\end{acknowledgements}

%



%
%

\end{document}